\documentclass[a4paper,11pt,reqno]{amsart}
\usepackage{fullpage}
\usepackage{amsmath}
\usepackage{amsthm}
\usepackage{amssymb}
\usepackage{physics}
\usepackage{graphicx}
\usepackage[toc,page]{appendix}
\usepackage{floatrow}
\usepackage{mathtools}
\usepackage{tikz-cd} 
\usepackage{todonotes}
\usepackage{bbm}
\usepackage{soul}
\usepackage[english]{babel}
\usepackage[T1]{fontenc}
\usepackage{todonotes}
\usepackage{framed}
\usepackage[all]{xy}
\allowdisplaybreaks

\usepackage{hyperref}
\makeatletter
\providecommand\@dotsep{5}
\def\listtodoname{List of Todos}
\def\listoftodos{\@starttoc{tdo}\listtodoname}
\makeatother

\newtheorem{theorem}{Theorem}[section]
\newtheorem{lemma}[theorem]{Lemma}
\newtheorem{prop}[theorem]{Proposition}
\newtheorem{cor}[theorem]{Corollary}

\newtheorem*{thm*}{Theorem}

\theoremstyle{definition}
\newtheorem{definition}[theorem]{Definition}
\newtheorem{example}[theorem]{Example}

\theoremstyle{remark}
\newtheorem{remark}[theorem]{Remark}

\newcommand{\N}{\mathbb{N}}

\newcommand{\C}{\mathbb{C}}

\newcommand{\Toe}{\mathbb{T}}
\newcommand{\Pim}{\mathbb{O}}

\newcommand{\Hilbf}{\mathrm{Hilb}_f}
\newcommand{\Hilf}{\mathrm{Hilb}^{1}_{f}}
\newcommand{\SPSf}{\mathrm{SPS}_{f}^{\mathbb{C}}}

\newcommand{\Fock}{\mathcal{F}}

\newcommand{\Cs}{C$^*$-}

\DeclarePairedDelimiter\inp{\langle}{\rangle}
\DeclareMathOperator\coker{coker}

\title[Quadratic subproduct systems]{Quadratic subproduct systems, free products, \\and their C*-algebras}

\date{\today}
\author{Francesca Arici}
\address{Mathematical Institute, Leiden University, P.O. Box 9512, 2300 RA
Leiden, the Netherlands}
\email{f.arici@math.leidenuniv.nl}
\thanks{Partially funded by the Netherlands Organisation of Scientific Research (NWO) under grants 016.Veni.192.237, and VI.Vidi.233.231. Part of the EU Staff Exchange project 101086394 “Operator Algebras That One Can See”.}

\author{Yufan Ge}
\address{Mathematical Institute, Leiden University, P.O. Box 9512, 2300 RA
Leiden, the Netherlands}
\curraddr{}
\email{y.ge@math.leidenuniv.nl}
\begin{document}

\begin{abstract}
Motivated by the interplay between quadratic algebras, noncommutative geometry, and operator theory, we introduce the notion of quadratic subproduct systems of Hilbert spaces. Specifically, we study the subproduct systems induced by a finite number of complex quadratic polynomials in noncommuting variables, and describe their Toeplitz and Cuntz--Pimsner algebras. Inspired by the theory of graded associative algebras, we define a free product operation in the category of subproduct systems and show that this corresponds to the reduced free product of the Toeplitz algebras. Finally, we obtain results about the K-theory of the Toeplitz and Cuntz--Pimsner algebras of a large class of quadratic subproduct systems.
\end{abstract}
    
\maketitle

\section{Introduction}
The study of subproduct systems and their C*-algebras has become a significant area of research at the intersection of multivariate operator theory \cite{Arv}, noncommutative geometry, and operator algebras. First introduced by Shalit and Solel in \cite{ShSo09}, and around the same time by Bhat and Mukherjee in the Hilbert space setting \cite{BM}, subproduct systems provide a natural framework for understanding row-contractive tuples of operators subject to polynomial constraints.

In this paper, we focus on \emph{quadratic} subproduct systems, which are subproduct systems of Hilbert spaces arising from a finite set of quadratic polynomials in a finite number of noncommuting variables. This class exhibits rich algebraic and operator-theoretic properties, and is quite a natural one to consider, given that algebras are often given in terms of commutation rules between their generators. Indeed, noncommutative algebras defined by quadratic relations are crucial examples of noncommutative spaces, such as those appearing in Manin's programme for noncommutative geometry \cite{Ma88,Ma91}. Such quadratic algebras include the deformations of quantum groups---and spaces---arising from an $R$-matrix, as defined in the seminal work of Faddeev, Reshetikhin, and Takhtajan \cite{FRT}. These continue to play a central role in noncommutative geometry, providing a rich source of examples of noncommutative spaces.

The interaction between subproduct systems and both classical and quantum groups extends beyond their mere construction, offering insights into the algebraic, geometric, and topological aspects of the underlying noncommutative spaces \cite{AK21,HaNe21,HaNe22,AGN24}. The presence of quantum group symmetries allows for elegant computations of the K-theoretic invariants of their C*-algebras. More recently, Aiello, Del Vecchio, and Rossi have introduced a subproduct system of finite-dimensional Hilbert spaces associated to the Motzkin planar algebra \cite{AdVR}, generalising the Temperley--Lieb subproduct systems of Habbestad and Neshveyev \cite{HaNe21,HaNe22}. 

Building on the framework established in previous studies, we examine the subproduct system analogue of the free product construction for noncommutative associative algebras. One of our motivations, in addition to the naturality of the free product construction, comes from the representation theory of the quantum group $SU_q(2)$. Our starting point is the observation that a free-product structure naturally appears when applying the algorithm in \cite[Section~2]{AK21} to multiplicity-free representations. This feature allows us to derive new insights into the algebraic and analytical properties of such subproduct systems, particularly in the context of their Fock spaces and associated \Cs algebras. 

The structure of the paper is as follows: We start by recalling the basic definitions and constructions for subproduct systems in Section~\ref{sec:sps}. In Section~\ref{sec:quadratic}, we introduce quadratic subproduct systems, highlighting their connections with quadratic algebras. We also define quadratic subproduct systems with few relations and discuss their relation to generic quadratic algebras, for which a lot is known about their growth and Hilbert series.   Section~\ref{sec:free} is devoted to the free product operation on subproduct systems. Here, we establish explicit formulas for the fibres of the free product and describe their fusion rules. Our main result for the section, Proposition~\ref{prop:fp_HS}, is a decomposition theorem for the Fock space of the free product of subproduct systems.

In Section~\ref{sec:Toe}, we study the Toeplitz algebras associated with these free products. Theorem~\ref{thm:toeplitz_free} asserts that the free product structure is preserved at the level of Toeplitz algebras, more precisely in terms of a reduced free product. Moreover, Theorem~\ref{thm:fp_kk-eq}, stated below, ensures that the free product construction allows us to bootstrap properties such as nuclearity and KK-equivalence to the complex numbers from smaller building blocks.
\begin{thm*}
    Let $\mathcal{H}$ and $\mathcal{K}$ be standard subproduct systems of Hilbert spaces. Assume that the Toeplitz algebras $\Toe_{\mathcal{H}}$ and $\Toe_{\mathcal{K}}$ are both \emph{nuclear} and KK-equivalent to the complex numbers. Then so is the Toeplitz algebra $\Toe_{\mathcal{H}\star \mathcal{K}}$.
\end{thm*}
We also demonstrate how our free product construction applied to monomial quadratic ideals corresponds to the graph join operation at the level of Cuntz--Pimsner algebras. 

Section~\ref{sec:TL} focuses on Temperley–Lieb subproduct systems, a subclass of quadratic subproduct systems defined by specific combinatorial constraints, introduced in \cite{HaNe21} and further studied in \cite{HaNe22,AGN24}. We analyse the free products of Temperley–Lieb subproduct systems, compute their K-theory, and construct explicit KK-equivalences for their Toeplitz algebras. We conclude the paper by studying the subproduct system of a finite-dimensional unitary $SU_q(2)$-representation, answering some questions regarding their structure and K-theory that were left open in \cite{AK21}.

\textbf{Acknowledgements.} We are pleased to thank Dimitris Gerontogiannis, Marcel de Jeu, Bram Mesland, Sergey Neshveyev, and Adam Rennie for interesting discussions. Jens Kaad deserves a special mention for having inspired and encouraged this research, and so does Wout Gevaert for his Master's thesis research involving computations we use in Section~\ref{sec:TL}. Finally, FA would like to thank Bernd Sturmfels for having made her aware of the beauty of quadratic relations\footnote{Never before has she been so sure that quadratic polynomials are extremely practical, and as such, good!}, and Tatiana Gateva-Ivanova for her introduction and guidance into the realm of noncommutative associative algebras.

\section{Preliminaries on subproduct systems and their algebras}
\label{sec:sps}
We start this section by recalling some basic facts from the theory of subproduct systems of Hilbert spaces and their \Cs algebras. Our main references are \cite{ShSo09, Vis12}. Although in their original paper, Shalit and Solel studied subproduct systems in the more general setting of \Cs and W$^*$-correspondences, we shall focus here on the Hilbert space case.

By a \emph{subproduct system} of finite-dimensional Hilbert spaces, we shall mean a sequence of finite dimensional Hilbert spaces $\mathcal{H} = \{H_n\}_{n \in \mathbb{N}_{0}}$, $\dim(H_0)=1$, together with isometries 
$$
\iota_{m,n}\colon H_{m+n} \to H_m \otimes H_n, 
$$ 
satisfying
\[ (\iota_{m,n}\otimes 1) \iota_{m+n,k} = (1 \otimes \iota_{n,k}) \iota_{m,n+k} \colon H_{m+n+k} \to H_m \otimes H_n \otimes H_k, \]
for all $m,n,k \in \mathbb{N}_{0}$, where $1$ denotes the identity operator. 

% \[\Pim_{\mathcal{H}} = \Toe_{\mathcal{H}}/\K(\mathcal{F}_\mathcal{H}).
% \]

A subproduct system is called \emph{standard} if $H_0=\C$, $H_{m+n}  \subseteq H_{m} \otimes H_{n}$, and the maps $
\iota_{m,n}$ agree with the embedding maps. 

As pointed out in \cite{ShSo09}, standard subproduct systems of finite-dimensional Hilbert spaces provide the natural framework for studying row-contractive tuples of operators subject to polynomial constraints, as made transparent by the existence of a \emph{noncommutative Nullstellensatz}.

\begin{prop}[{\cite[Proposition 7.2]{ShSo09}}]\label{prop:NSS}
Let $H$ be a $d$-dimensional Hilbert space. Then there is a bijective inclusion-reversing correspondence between the proper homogeneous ideals $J  \subset \mathbb{C} \langle X_1,\ldots, X_d \rangle$ and the standard subproduct systems $\lbrace H_n \rbrace_{n \in \mathbb{N}_{0}}$ with $H_1 \subseteq H$.
%(all structure maps are given by canonical inclusions).
\end{prop}
Let us fix an orthonormal basis $\lbrace e_i \rbrace_{i=1}^d$ for $H$. The correspondence works as follows: for a noncommutative polynomial $P=\sum c_{\alpha} X^{\alpha}$ in variables $X_1\dots,X_d$, we write $P(e) = \sum c_{\alpha} e_{\alpha}$, where $e_{\alpha}=e_{\alpha_1} \otimes \dots \otimes e_{\alpha_k}$ for $\alpha = \alpha_1 \cdots \alpha_k$ a length-$k$ word. To any proper homogeneous ideal $J \subset \C \langle X_1,\ldots,X_d \rangle$, one associates the standard subproduct system with fibres $H_n:= H^{\otimes n} \ominus \lbrace P(e) \vert P \in J^{(n)} \rbrace$, for every $n \geq 0$, where $J^{(n)}$ denotes the degree-$n$ component of $J$. 

Following \cite[Definition 7.3]{ShSo09}, we refer to $\mathcal{H}^J$ and $J_{\mathcal{H}}$ as the \emph{subproduct system associated with the ideal $J$} and the \emph{ideal associated with the subproduct system $\mathcal{H}$}, respectively.

While, in principle, the above construction depends on the choice of an orthonormal basis for $H$, different choices yield isomorphic subproduct systems in the sense of \cite[Definition~1.4]{ShSo09}.

\begin{prop}[{\cite[Proposition~7.4]{ShSo09}}]
\label{prop:unit_iso}
     Let $\mathcal{H}$ and $\mathcal{K}$ be standard subproduct systems with $\dim(H_1) = \dim(K_1) = d <\infty $. Then $\mathcal{H}$ is isomorphic to $\mathcal{K}$ if and only if there is a unitary linear change of variables in $\mathbb{C}\langle X_{1},\ldots, X_{d}\rangle$ that sends $J_{\mathcal{H}}$ onto $J_{\mathcal{K}}$. 
\end{prop}

In a basis-independent fashion, Proposition~\ref{prop:NSS} can also be formulated as follows: There is a bijective inclusion-reversing correspondence between the proper homogeneous ideals $J$ of the free algebra in $\dim(H)$-generators and the standard subproduct systems $\lbrace H_n \rbrace_{n \in \mathbb{N}_{0}}$ with $H_1 \subseteq H$.

It is worth recalling that all standard subproduct systems are obtained this way, see \cite[Proposition 7.2]{ShSo09}.  As we mentioned in the introduction, we will focus on standard subproduct systems induced by a finite number of quadratic polynomials in non-commuting variables, as these form a more tractable class of examples.

\subsection{Toeplitz and Cuntz--Pimsner algebras of subproduct systems}
We conclude this section by recalling the construction of the Toeplitz and Cuntz--Pimsner algebras of a subproduct system of Hilbert spaces. 

The Fock space of the subproduct system $\mathcal{H}$ is the direct sum Hilbert space 
\[\mathcal{F}_{\mathcal{H}} := \bigoplus_{n\ge0} H_n.\]
On  the Hilbert space $\mathcal{F}_{\mathcal{H}}$ we consider operators defined by
\[ T_\xi(\zeta) := \iota_{1,n}^*(\xi\otimes \zeta), \quad \xi \in H_1,\ \zeta \in H_n.\]
Note that the Fock space is a subspace of the full Fock space of $H_1$:
\[ \mathcal{F}_{\mathcal{H}} \subseteq \bigoplus_{n\ge0} H_1^{\otimes n} \quad \mbox{ and } \quad T_\xi(\zeta) = f_{n+1}(\xi\otimes \zeta), \quad \xi \in H_1, \zeta \in H_n, \]
where $f_{n+1}$ is the projection $f_{n+1} \colon H_1^{\otimes (n+1)} \to H_{n+1}$.

The \emph{Toeplitz algebra} $\Toe_{\mathcal{H}}$ of the subproduct system $\mathcal{H}$ is the unital \Cs algebra generated by $T_1,T_2,...,T_d$, where $T_i = T_{\xi_i}$ for an orthonormal basis $(\xi_i)_{i=1}^d$ of $H_1$. If one denotes by $e_0$ the rank-one projection onto $H_0$, it is straightforward to verify that \[ 1_{\mathcal{F}} - \sum_{i=1}^{d} T_iT_i^* = e_0.\] 
Consequently, the compact operators on the Fock space $\mathbb{K}(\mathcal{F}(\mathcal{H}))$ are contained in~$\Toe_{\mathcal{H}}$ (cf. \cite[Corollary 3.2]{Vis12}). This fact is used to define the \emph{Cuntz--Pimsner algebra} $\Pim_{\mathcal{H}}$ of the subproduct system as the quotient:
\begin{equation}\label{eq:CPext}
\xymatrix{0 \ar[r] & \mathbb{K}(\mathcal{F}_{\mathcal{H}}) \ar[r] & \Toe_{\mathcal{H}}\ar[r] & \Pim_{\mathcal{H}} \ar[r] & 0}.
\end{equation}%

\section{Quadratic subproduct systems from quadratic algebras}
\label{sec:quadratic}
\subsection{Quadratic Algebras and Their Hilbert Series}
In this work, we shall use several results from the theory of quadratic algebras, particularly in connection with their Hilbert series. Our main references are \cite{PoPo,Ufr}. Our base field will be the complex numbers.

Given a vector space $V$, we denote its tensor algebra by $\mathcal{T}(V)$. This is naturally graded by rank, and we write \[ \mathcal{T}(V):=\bigoplus_{n\geq 0} \mathcal{T}^n(V), \qquad \mathcal{T}^n(V):=V^{\otimes n}.\]

\begin{definition}
A graded algebra $A \simeq \bigoplus_{n \geq 0} A_n$ is called \emph{one-generated} if the natural map $p: \mathcal{T}(A_1) \to A$ from the tensor algebra generated by $A_1$ is surjective. We call a one-generated algebra \emph{quadratic} if the ideal $J_A := \ker(p)$ is generated, as a two-sided ideal, by 
\[ I_A := J_A \cap \mathcal{T}^2(A_1) \subseteq A_1 \otimes A_1. \]
\end{definition}

In other words, a quadratic algebra is determined by a vector space of generators $V = A_1$ and a subset of relations $ I_A \subset A_1 \otimes A_1.$ We shall denote with $R$ the complex vector subspace of $A_1 \otimes A_1$ spanned by the relations. If $\dim(R) = r$, we call $A$ an \emph{$r$-relator} quadratic algebra.

Recall that for a graded vector space $V = \bigoplus_{k \geq 0} V_k$, with finite-dimensional graded components, its Hilbert series is the formal power series
\begin{equation}
\label{eq:Hilb_series}
h_V(z) = \sum_{k \geq 0} \dim (V_k) z^k.
\end{equation}

Hilbert series of associative algebras provide information about their growth. In \cite{DA83}, Anick studied them under certain finiteness hypotheses by considering a well-ordering defined as follows. Given two formal power series $f(z), g(z) \in \mathbb{R} [\![ z ]\!]$, we write $f(z) \geq g(z)$ if the inequality holds coefficient-wise. Moreover, we write $\vert f(z) \vert$ for the series obtained by deleting all the terms starting from the first negative term. Using this notation, one can write a lower bound for the Hilbert series of a quadratic algebra:

\begin{prop}{\cite[{Proposition~2.3}]{An88}}
For any quadratic algebra in $m$ generators and $r$ relations, the Hilbert series satisfies
\[ h_A(z) \geq \vert (1 - mz + rz^2)^{-1} \vert. \]
\end{prop}

\begin{theorem}[{cf. \cite[Proposition~4.1]{PoPo}}]
The minimal possible value of $\dim(A_3)$ for a quadratic algebra with $\dim(A_1) = d$ and $\dim(A_2) = s$ is
\[
\begin{cases}
  0 & \text{if } s \leq \frac{d^2}{2}, \\
  2ds - d^3 & \text{if } s \geq \frac{d^2}{2}.
\end{cases}
\]
\end{theorem}
Anick answered the question of which algebras attain the minimal Hilbert series by considering the question of genericity. Let us first clarify what we mean by the term \emph{generic}.

\begin{definition}[{cf. \cite[Chapter~6]{PoPo}}]
A complex \emph{generic} quadratic algebra in $d$ generators and $r$ relations is a generic point in the variety of quadratic algebras $\mathcal{Q}_{d,r}$ with $\dim(A_1) = d$ and $\dim(R) = r$.
\end{definition}

By definition, such an $r$-relator quadratic algebra is determined by a $(d^2 - r)$-dimensional subspace of $\mathbb{C}^{d^2}$. We can therefore identify the variety $\mathcal{Q}_{d,r}$ of such quadratic algebras with the complex Grassmannian $\mathbf{Gr}_{\mathbb{C}}(r, d^2).$. One then says that a \emph{generic} complex quadratic algebra is a quadratic algebra corresponding to a generic point in the Grassmannian variety of quadratic algebras (cf. \cite[Lemma~4.1]{DA83}).

Anick's main result establishes that generic algebras are exactly those which possess the coefficient-wise minimal Hilbert series \cite[Definition~4.9]{DA83}.
Additionally, in the quadratic case, adding some further constraints on the number of generators and relations allows one to obtain an explicit formula for the generic Hilbert series:

\begin{prop}[{\cite[Proposition~4.2]{PoPo}}]
\label{prop:HilbSeries}
A generic quadratic algebra $A$ in $d$ generators and $r$ relations is Koszul\footnote{Koszulness, in algebra, refers to a property of a graded or filtered algebra where its minimal free resolution has certain desired homological properties. We shall not go into the precise definition here, but refer the reader to \cite{PoPo} for more details.} if and only if one of the inequalities holds:
\[ r \geq \frac{3 d^2}{4}, \qquad r \leq \frac{d^2}{4}. \]
Then the Hilbert series of $A$ is either 
\begin{align}
\label{eq:HSgenFin}
h_A(z) &= 1 + dz + (d^2 - r)z^2, & \text{or} \\
\label{eq:HSGeneric}
h_A(z) &= (1 - dz + rz^2)^{-1}, &
\end{align}
respectively.
\end{prop}
Note that in the first case, when $r \geq {3 d^2}/{4}$, \eqref{eq:HSgenFin} implies that the quadratic algebra is finite-dimensional (and in particular $\dim(A_3)=0$). This should serve as a motivation for focusing on the constraint $r \leq {d^2}/{4}$ later on. In that case, we shall talk about a quadratic algebra with \emph{few relations}. Let us stress that, by Anick's Theorem, a quadratic algebra $A$ with few relations is generic if and only if its Hilbert series equals \eqref{eq:HSGeneric}.

\begin{remark}
In the literature, one can encounter another related notion of having few relations: \cite{Zh97}. The main result there states that whenever $\mathrm{rank}(R) > \dim(R) + 1$, the quadratic algebra is Koszul of global dimension $2$, and its Hilbert series is given by \eqref{eq:HSGeneric}, making such algebras automatically generic.
\end{remark}

Hilbert series of quadratic algebras with a fixed number of generators and relations are well-studied and understood, at least in low dimensions. We refer the reader to \cite[Section~6.5]{PoPo} for some explicit expressions of the Hilbert series.
\begin{example}
When $d = 2$ and $r = 1$, a \emph{non-generic} quadratic algebra must necessarily have the following Hilbert series:
\begin{equation}
\label{eq:HSnonGen}
(1 - z - z^2)^{-1} = 1 + 2z + 3z^2 + 5z^3 + 8z^4 + 13z^5 + \dots.
\end{equation}
An example is the algebra $\mathbb{C}\langle X_1, X_2 \rangle / \langle X_1^2 \rangle$, which we will re-encounter in Example~\ref{ex:fibonacci}.
\end{example}

\subsubsection{Free products of quadratic algebras and their Hilbert series}

The category of unital graded algebras over a field has a natural coproduct operation, given by the algebraic free product. 

\begin{definition}
    Given two graded algebras $A$ and $B$ over a field ${k}$, their free product, denoted $A \sqcup B$, is defined as the associative algebra generated freely by $A$ and $B$. Explicitly:
\begin{equation}
\label{eq:free_prod_alg}
    A \sqcup B := \bigoplus_{i\geq 0, \epsilon_1, \epsilon_2 \in \lbrace 0, 1 \rbrace} A_{+}^{\otimes_{\epsilon_1}}\otimes (B_{+} \otimes A_{+})^{\otimes i} \otimes B_{+}^{\otimes_{\epsilon_2}},
\end{equation}
with the usual convention that $A_{+}:=\bigoplus_{n \geq 1} A_n$.
\end{definition}
It is natural to wonder how the Hilbert series behaves when one considers free products. Given two finitely-presented algebras $A$ and $B$, by \cite[Theorem~4.5.3]{Ufr}, the Hilbert series of the free product algebra $A\sqcup B$ can be expressed in terms of the Hilbert series of $h_A(z)$ and $h_B(z)$ of the algebras $A$ and $B$:
\begin{equation}
\label{eq:HSprod}
    h_{A\sqcup B}(z)^{-1} = h_A(z)^{-1} + h_B(z)^{-1} +1 .
\end{equation}
This has important consequences for the question of genericity.
\begin{remark}
\label{rem:genericProd}
    Suppose that $A$ and $B $ are $r_0$- and $r_1$-relator generic quadratic algebras in $d_0$ and $d_1$ generators, respectively, satisfying the additional condition $r_i \leq d_i^ 2/ 4,$ for $i=0,1$.  
% \begin{equation}
% \label{eq:gfp}
%     (E^{(0)} \star E^{(1)})_m \otimes (E^{(0)} \star E^{(1)})_1= (E^{(0)} \star E^{(1)})_{m+1} \oplus (E^{(0)} \star E^{(1)})_{m-1}^{\oplus (R_0 + R_1)} 
% \end{equation}
The formula for the Hilbert series of the free product of algebras \eqref{eq:HSprod} yields \[ h_{A \sqcup B} (z) = (1 - (d_0 + d_1)z + (r_0 + r_1)z^2)^{-1},\]
    implying that $A\sqcup B$ is a generic quadratic algebra in $d_0+d_1$ generators and $r_0 + r_1$ relations.
\end{remark}

Other important operations and constructions that preserve the class of quadratic algebras are Veronese powers and Segre products \cite{PoPo}. We defer the treatment of their operator algebraic counterparts to future work.

\subsection{Quadratic subproduct systems}

Having discussed the fundamentals of the theory of quadratic algebras, we are ready to introduce quadratic subproduct systems of Hilbert spaces.

Recall first the definition of the maximal suproduct system with prescribed fibres.

\begin{definition}[{\cite[Section~6]{ShSo09}}]
\label{defn:max}
Let $H_0 = \C$, and let $H_1$ be a Hilbert space. Let $H_i$,$i=2, \dots, k$ be subspaces of $H_1^{\otimes i}$ respectively. The maximal subproduct system with prescribed fibres up to level $k$ is defined as the subproduct system $\mathcal{H}=\lbrace H_n \rbrace_{n\geq 0}$ with   
\begin{align*}
	H_n = \bigcap_{i + j = n} H_i \otimes H_j, \qquad n >k.
\end{align*}
\end{definition}

\begin{prop}
\label{prop:quad}
Let $\mathcal{H}$ be a standard subproduct system of finite-dimensional Hilbert spaces with $\dim(H_1)=d$. The following facts are equivalent: 
\begin{enumerate}
    \item{\label{max}} There exists $R \subseteq H_1 \otimes H_1$ such that $H$ is isomorphic to the \emph{maximal} subproduct system with fibres $\C, H_1, R^\perp$ in the sense of Definition~\ref{defn:max}.
    \item There exist finitely many homogeneous quadratic polynomials $f_1, \dots, f_k \in \C \langle X_1, \dots, X_d \rangle$ such that $\mathcal{H} \simeq H^{\langle f_1,\dots,f_k \rangle}$.
    \item The ideal $J_{\mathcal{H}}$ is generated by quadratic polynomials. 
\end{enumerate}
\end{prop}
A standard subproduct system of finite-dimensional Hilbert spaces satisfying any of the above conditions will be called a \emph{quadratic} subproduct system of Hilbert spaces.
\begin{definition}
    A quadratic subproduct system with $\dim(H_1) = d$ and $\dim(R)=r$ will be called an \emph{$r$-relator} quadratic subproduct system.  If $R$ corresponds to a generic point of the Grassmannian $\mathbf{Gr}_{\mathbb{C}}(r,d^2)$, then we call $\mathcal{H}$ a \emph{generic} $r$-relator quadratic subproduct system. 
\end{definition}
Note that everything still makes sense in the case $r=0$, as one obtains product systems of finite-dimensional Hilbert spaces. Everything we discuss in the rest of this work holds in that setting as well. 

A special class of examples is that of \emph{one-relator} quadratic subproduct systems, i.e., those subproduct systems whose underlying ideal $J_{\mathcal{H}}$ is generated by a single quadratic polynomial. In that case, we have the following result due to Shalit and Solel.
\begin{theorem}[{\cite[Proposition 11.1]{ShSo09}}]
 Let $A,B \in M_d(\mathbb{C})$. Consider the two quadratic subproduct systems $\mathcal{H}_A$ and $\mathcal{H}_{B}$ given by the polynomials 
 \[\sum_{i,j=1}^d A_{ij}X_iX_j, \qquad \sum_{i,j=1}^d B_{ij}X_iX_j, \]
respectively. Then there is an isomorphism $V: \mathcal{H}_A \to \mathcal{H}_B$ if and only if there exists $\lambda \in \mathbb{C}$ and a unitary $d \times d$ matrix $U$ such that $B = \lambda U^{t}AU$. 
\end{theorem}

In the following, when talking about the Hilbert series of a subproduct system of Hilbert spaces, we will mean the formal power series \eqref{eq:Hilb_series} of the underlying graded vector space.

Motivated by \ref{prop:HilbSeries}, we give the following definition.

\begin{definition}
\label{def:gen_few}
Let $\mathcal{H}$ be an $r$-relator quadratic subproduct system of finite-dimensional Hilbert spaces with $\dim(H_1)=d$.  We say that $\mathcal{H}$ has \emph{few relations} if 
\begin{equation}
\label{eq:few}
r \leq \frac{d^2}{4}.
\end{equation}
\end{definition}

\begin{remark}
\label{rem:intersection}
For any quadratic subproduct system of finite-dimensional Hilbert spaces $\mathcal{H}$,
\begin{equation}
\label{eq:def_fibres}
    H_{m+1} = H_{1}\otimes H_{m} \cap H_{m}\otimes H_{1}, \qquad \mbox{for all }  m \geq 2.
\end{equation}
\end{remark}

Indeed, a quadratic subproduct system is a maximal standard subproduct system with prescribed fibres $H_1$ and $H_2$, for which the higher-level fibres are given as in Definition~\ref{defn:max}.
Applying this to $H_m\otimes H_1$ and $H_1 \otimes H_m$, we obtain 
\begin{align*}
    H_m\otimes H_1 = \bigcap_{k + l = m} H_k \otimes H_l \otimes H_1, \\
    H_1 \otimes H_m = \bigcap_{k + l = m} H_1 \otimes H_k \otimes H_l,
\end{align*}
from which we obtain that $H_i \otimes H_{m+1 - j} = H_i \otimes H_{m-i} \otimes H_1 \cap H_1 \otimes H_{i-1} \otimes H_{m+1 - i}$, since $H_i \subset H_1 \otimes H_{i-1}$ and $H_{m+1-i} \subset H_{m-i}\otimes H_1$. By simple linear algebra arguments, \eqref{eq:def_fibres} follows.

\begin{remark}One may wonder whether the definition of a quadratic subproduct system could be extended to the setting of correspondences. This is clearly the case if we use Condition~\ref{max} in Proposition~\ref{prop:quad} and consider maximal subproduct systems with prescribed fibres in degrees up to two. However, two important aspects require additional care. Firstly, in the correspondence case, one loses the connection to the theory of polynomials in non-commuting variables. Most importantly, there are some technical issues related to the complementability of submodules, deeply connected to the theory of two-projections in Hilbert modules \cite{MeRe}. We postpone the treatment of quadratic subproduct systems of correspondences to future work.
\end{remark}

The subproduct systems associated with an irreducible $SU(2)$-representation from \cite{AK21}, and more generally the Temperley--Lieb subproduct systems studied in \cite{HaNe21, HaNe22, AGN24} are clear examples of one-relator quadratic subproduct systems, being defined by a single quadratic polynomial.

\subsubsection{The importance of being generic}
As anticipated, we shall focus on \emph{generic} subproduct systems with few relations, that is, satisfying
$ r \leq {d^2}/{4}$.

\begin{prop}
Let $\mathcal{H}$ be a generic quadratic subproduct system that satisfies the assumptions in Definition~\ref{def:gen_few}. 
Denote $\delta_n = \dim(H_n)$.  The sequence $\lbrace \delta_n \rbrace_{n \geq 0}$ satisfies the following recurrence relation:
    \begin{equation}
        \delta_0  =1;\quad \delta_1  =d; \quad \delta_{n+1}  = d\cdot \delta_n - r \cdot \delta_{n-1}, \qquad n \geq 1.
    \end{equation}
    \end{prop}
\begin{proof}
The proof relies on a standard argument that involves the logarithmic derivative of the generating function for the Hilbert series of $\mathcal{H}$. 
\end{proof}

As a consequence, we have the following:

\begin{cor}
Let $\mathcal{H}$ be a generic quadratic subproduct system that satisfies the assumptions in Definition~\ref{def:gen_few}. Then, for every $n \geq 1$, there are vector space isomorphisms
\[H_n \otimes H_1  \cong {H}_{n+ 1} \oplus H_{n - 1}^{\oplus r}, \]
\end{cor}

Finding explicit formulas for isometries that implement the above isomorphism may be a hard task. This is possible in some cases, including those we will encounter in Sections \ref{sec:TL} and \ref{sec:SUq2}.
\begin{example}[The Fibonacci subproduct system]
\label{ex:fibonacci}
Consider the subproduct system associated to the ideal $I:= \langle X_1^2 \rangle \subseteq \mathbb{C}\langle X_1, X_2 \rangle$. It is easy to see that it has Hilbert series \eqref{eq:HSnonGen}, and as such, it cannot be generic. The dimension sequence of this subproduct system is the (shifted) Fibonacci sequence: 
\[ \delta_0=1, \quad \delta_1=2, \quad \delta_{n+1}= \delta_{n} + \delta_{n-1}, \qquad n \geq 1.\]
\end{example}
\section{Free products of subproduct systems and their Fock spaces}
\label{sec:free}
We are ready to define the free product of quadratic subproduct systems:
\begin{definition}
\label{defn:free_quad_SPS}
	Let $\mathcal{H}= \{ H_n \}_{n \in \mathbb{N}_{0}} $ and $\mathcal{K}= \{K_n \}_{n \in \mathbb{N}_{0}}$ be two quadratic subproduct systems.  We then define the free product $\mathcal{H} \star \mathcal{K}$
 of $\mathcal{H}$ and $\mathcal{K}$ 
 as the maximal subproduct system with prescribed fibres $(H\star K)_1 := H_1 \oplus K_1, (H\star K)_2 := (H_2^\perp \oplus K_2^\perp)^\perp$.
\end{definition}
In the polynomial picture, if we write  $J_{\mathcal{H}}$ and $J_{\mathcal{K}}$ for the corresponding ideals, then the free product is obtained by considering the ideal generated by the disjoint union of the quadratic polynomials that generate $J_{\mathcal{K}}$ and $J_{\mathcal{K}}$ in the free algebra in $\dim({H_1})+\dim(K_1)$ variables.

We shall use Remark~\ref{rem:intersection} to describe the Hilbert space $(H \star K)_2$ explicitly.

\begin{lemma}\label{lem:fibre2}
	Let $\mathcal{H}$ and $\mathcal{K}$ be two quadratic subproduct systems of finite-dimensional Hilbert spaces. We have an isomorphism of inner product spaces \[ (H \star K)_2 \simeq
 H_2 \oplus K_2 \oplus (H_1\otimes K_1) \oplus (K_1 \otimes H_1).\]
\end{lemma}
We are interested in describing the fibres of the free product of quadratic subproduct systems using the formula from Remark~\ref{rem:intersection}, which involves tensor products and intersections. Let us first recall some enumerative combinatorics.
\begin{definition}
    A \emph{composition} of a positive integer $n \geq 1$ is a sequence $\sigma_1, \dots, \sigma_p $, $\sigma_i \in \mathbb{Z}_{+}$ with $\sum_{i=1}^p \sigma_i = n$. The $\sigma_i$'s are called the \emph{parts} of $n$. We denote the set of compositions of $n$ with $\mathcal{C}(n)$ and the subset of compositions of $n$ into \emph{exactly} $p$ parts by $\mathcal{C}_p(n)$. \end{definition}
A given integer $n\geq 1$ has $2^{n-1}$ compositions. Moreover, for every $1 \leq p \leq n$ the cardinality of $\mathcal{C}_p(n)$ equals the binomial coefficient $ \binom{n-1}{p-1}$.

Let $\mathcal{H}^{(0)}$ and $\mathcal{H}^{(1)}$ be two subproduct systems of finite-dimensional Hilbert spaces, quadratic or not.  For a fixed $m \in \mathbb{N}_{0}$, $1 \leq p \leq m$ and $\underline{d} \in \mathcal{C}_p(m)$ and $j=0,1$, we define 
\[
H_{\underline{d}}^{(j)} :=   H^{(j)}_{d_1} \otimes H^{((j+1) \bmod{2})}_{d_2}  \otimes \cdots \otimes H^{((j+p +1) \bmod 2)}_{d_p}.\]

As already mentioned, we are interested in understanding how such spaces behave for tensor products and intersections, since those are the operations involved in the construction of a quadratic subproduct system, and in general, in the definition of a subproduct system with prescribed fibres. 

\begin{prop}
\label{prop:intersection}
Let $\underline{d} \in C_p(m)$, $\underline{f} \in C_q(m)$. With the above notation, the intersection \[ (H_1^{(i)} \otimes H^{(j)}_{\underline{d}} )\cap  (H^{(l)}_{\underline{f}} \otimes H_{1}^{(k)} )= \emptyset,  \]
whenever $i \neq l$, $q \neq p, p-1, p+1$.

Moreover, we have 
\[
( H_1^{(i)} \otimes H^{(j)}_{\underline{d}})\cap (H^{(i)}_{\underline{f}} \otimes H_{1}^{(k)} )= 
\begin{cases} 
H^{(i)}_{(d_1+1,d_2,\dots,1)}               & i=j, q=p-1, k = p+i+1 \bmod 2, d_p=1, \\ 
                                            & f_1=d_1+1, f_2=d_2, \dots f_{p-1}=d_{p-1}; \\
H^{(i)}_{(d_1+1,d_2,\dots,d_p)}            & i=j, q=p,  k = p+i+1 \bmod 2, d_p>1, \\
                                            & f_1=d_1+1, f_2=d_2, \dots f_p=d_p-1; \\
H^{(i)}_{(1,d_1,d_2,\dots,d_{p-1},1)}       & i \neq j, q=p, d_p=1, \\
                                            & f_1 =1, f_2=d_1, \dots, f_p=d_{p-1};\\  
H^{(i)}_{(1,d_1,d_2,\dots,d_{p-1}, d_p)}    & i \neq j , q=p+1, d_p > 1\\
                                            & f_1 =1, f_2=d_1, \dots, f_p=d_{p-1}, f_{p+1}=d_p+1;\\
\emptyset & \mbox{otherwise}
\end{cases}
\]

\end{prop}
\begin{proof}
    Let us first suppose that 
$i=j$, then we are looking at the subspace 
\[
H_1^{(i)} \otimes H^{(i)}_{\underline{d}} \cap H^{(i)}_{\underline{f}} \otimes H_{1}^{(k)},
\]
which we can rewrite as 
\[H_1^{(i)} \otimes H^{(i)}_{d_1} \otimes H^{((i+1)\bmod{2})}_{d_2} \otimes \cdots \otimes H^{((p+i+1)\bmod{2})}_{d_p} \cap H^{(i)}_{f_1} \otimes H^{((i+1)\bmod 2)}_{f_2} \otimes \cdots H_{f_q}^{((q+i+1)\bmod 2)} \otimes H_{1}^{(k)}.\]

First of all, we observe that if $p=1$, the intersection amounts to 
\[
 H_1^{(i)} \otimes H^{(i)}_{d_1} \cap H^{(i)}_{f_1} \otimes H_{1}^{(k)} = \begin{cases}
 H^{(i)}_{d+1} & f=d>1,  k=j=i, \\
 \emptyset & \mbox{otherwise.}
 \end{cases}
 \]

If $p>1$, the non-triviality of the intersection forces $\underline{d}$ and $\underline{f}$ to satisfy either \[
d_1+1=f_1, d_2=f_2, \dots d_{p-1}=f_q, d_p=1,\] 
for $q=p-1$, $k = (p+i+1) \bmod 2$; or
\[
d_1+1=f_1, d_2=f_2, \dots d_{p-1}=f_{q-1}, d_p=f_q+1>1,\]
for $p=q$ and $k = (p+i+1) \bmod 2$.

Similar considerations give the claim for the intersections for $i \neq j$.
\end{proof}
With Proposition~\ref{prop:intersection} in place, we can now provide an explicit description of the fibres of the free product of two quadratic subproduct systems.

\begin{prop}
\label{prop:fibres_product}
Let $\mathcal{H}^{(0)}$ and $\mathcal{H}^{(1)}$ be two quadratic subproduct systems. Their free product $\mathcal{H}^{(0)} \star \mathcal{H}^{(1)}$ satisfies 
\begin{equation}
\label{eq:decomp_free_prod}
    (H^{(0)} \star H^{(1)})_m= \bigoplus_{i\in \lbrace 0,1\rbrace} \bigoplus_{p=1}^{m} \bigoplus_{\underline{d} \in \mathcal{C}_p(m)}  H^{(i)}_{\underline{d}} .
\end{equation}    \end{prop}

Taking the free product of two quadratic subproduct systems is an associative operation:
\begin{lemma}
    Let $\{ H^{(i)}_m\}_{m \in \mathbb{N}_{0}}$, $i=0,1,2$ be three quadratic subproduct systems of finite-dimensional Hilbert spaces. Then for all $m \geq 0$, we have unitary isomorphisms \[((H^{(1)}\star H^{(2)})\star H^{(3)}))_m = (H^{(1)}\star (H^{(2)}\star H^{(3)}))_m .\]
\end{lemma}
Consequently, one can unambiguously consider the free product of a finite number of quadratic subproduct systems. 

\subsection{Fusion Rules for free product subproduct systems}
    Let $\mathcal{H}= \{ H_m\}_{m \in \mathbb{N}_{0}}$ and $\mathcal{K}=\{ K_m\}_{m \in \mathbb{N}_{0}} $ be $r_0$-relator and $r_1$-relator generic quadratic subproduct systems in $d_0$ and $d_1$ generators, respectively, satisfying the condition in \eqref{eq:few}. Then we can apply the same argument from Remark~\ref{rem:genericProd} to deduce that their free product subproduct system also satisfies Condition \eqref{eq:few} of having few relations, has Hilbert series \[ h_{\mathcal{H} \star \mathcal{K}} (z) = (1 - (d_0 + d_1)z + (r_0 + r_1)z^2)^{-1}.\]
Setting $\delta_m = \dim((H\star K)_m)$,  we obtain the sequence $\lbrace \delta_m \rbrace_{m \geq 0}$ satisfies the recurrence relation
    \begin{align*}
        \delta_0 & =1;\\
        \delta_1 & =(d_0+d_1);\\   
        \delta_{m+1} & = (d_0 + d_1)\cdot \delta_m - (r_0 + r_1)\cdot \delta_{m-1},
    \end{align*}
    which implies    
\begin{equation}
\label{eq:gfp}
    (H \star K)_m \otimes (H \star K)_1= (H \star K)_{m+1} \oplus (H \star K)_{m-1}^{\oplus (r_0 + r_1)} .
\end{equation}
% \begin{equation}
% \label{eq:gfp}
%     (E^{(0)} \star E^{(1)})_m \otimes (E^{(0)} \star E^{(1)})_1= (E^{(0)} \star E^{(1)})_{m+1} \oplus (E^{(0)} \star E^{(1)})_{m-1}^{\oplus (R_0 + R_1)} 
% \end{equation}

%     Note that $\mathcal{H}$ and $\mathcal{K}$ have the Hilbert series given by 
%         \begin{align*}
%         h_0 (t) = (1 - n_0 t + r_0t^2)^{-1}, \\
%         h_1 (t) = (1 - n_1 t + r_1 t^2)^{-1}, 
%     \end{align*}
%     Applying the formula for the Hilbert series of the free product of algebras \eqref{eq:HSprod},
%     we obtain 

% s the desired isomorphism.

If additionally  $\mathcal{H}$ and $\mathcal{K} $ are \emph{generic} one-relator quadratic subproduct systems, the isomorphism \eqref{eq:gfp} reduces to
\begin{equation}
\label{eq:freeprod}
    (H \star K)_m \otimes (H \star K)_1= (H \star K)_{m+1} \oplus (H \star K)_{m-1}^{\oplus 2} 
\end{equation}   

\subsection{Fock spaces of free products are free products of Fock spaces}
We start this subsection by discussing free products of Hilbert spaces. Our main references are \cite{EG2,Ger:Crelle,Sp98}.

In general, one can define the free product of a family of Hilbert spaces, but for the sake of readability, we shall focus here on the case of two Hilbert spaces only. 
\begin{definition}

\label{def:fpHilb}
Let $(H_1,\xi_0)$ and $(H_2,\xi_0)$ be two Hilbert spaces with a distinguished normal vector $\xi_0$. Their free product is the space $(H, \xi_0)$ with
	\begin{align*}
		H \simeq \mathbb{C} \xi_0 \oplus \bigoplus_{p \geq 1} \bigoplus_{i \in D_p} ( H_{i_1}^\circ \otimes \dots \otimes H_{i_p}^\circ),
	\end{align*}
	where 
	\begin{align*}
		D_p = \{ i = (i_1, i_2, \dots, i_p): i_j \in \{1, 2\} \text{ and } i_j \neq i_{j+1},\ \forall 1 \leq j \leq p-1 \},
	\end{align*}
and $H_{i}^\circ$ is the orthocomplement of $\mathbb{C} \xi_0$ in $H_{i}$.
\end{definition}
Applying Definition~\ref{def:fpHilb} to the Fock spaces of two quadratic subproduct systems, $\Fock_i := \Fock_{\mathcal{H}^{(i)}}$, with distinguished normal vector the vacuum vector $\omega_0$ for $H^{(1)}_0 \simeq \mathbb{C} \simeq H^{(2)}_0 $, we obtain 
 \begin{align*}
		\Fock_1 * \Fock_2 \simeq \mathbb{C} \omega_0 \oplus \bigoplus_{p \geq 1} \bigoplus_{i \in D_p} ( \mathcal{F}_{i_1}^+ \otimes \dots \otimes \mathcal{F}_{i_p}^+),
	\end{align*}
 with the usual convention that $\mathcal{F}_{i}^+$ denotes the positive Fock space, i.e. $\mathcal{F}(\mathcal{H}^{(i)})^+ = \bigoplus_{n\geq 1} H^{(i)}_n$.

\begin{prop}
\label{prop:fp_HS}
Let $\mathcal{H}$ and $\mathcal{K}$ be quadratic subproduct systems of finite-dimensional Hilbert spaces. Then the Fock space $\mathcal{F}(\mathcal{H} \star \mathcal{K})$ is unitarily isomorphic to the Hilbert space free product of the Fock spaces $\mathcal{F}({\mathcal{H}}) \ast \mathcal{F}(\mathcal{K})$.    
\end{prop}
\begin{proof}
First, we consider the free product of the Fock spaces 
	\[
		\mathcal{F}(\mathcal{H}) \ast \mathcal{F}(\mathcal{K})  = \mathbb{C} \oplus \bigoplus_{n \geq 1} \left( \underbrace{F(\mathcal{H})^+\otimes \mathcal{F}(\mathcal{K})^+\otimes \mathcal{F}(\mathcal{H})^+ \otimes \dots}_n \oplus \underbrace{\mathcal{F}(\mathcal{K})^+\otimes \mathcal{F}(\mathcal{H})^+\otimes \mathcal{F}(\mathcal{K})^+ \otimes \dots}_n \right).
\]
Any direct summand of $\mathcal{F}(\mathcal{H}) \ast \mathcal{F}(\mathcal{K})$ has the following form:
	\begin{equation}
    \label{eq:component}
		 H_{i_1} \otimes K_{j_1} \otimes H_{i_2} \otimes \dots \otimes H_{i_l} \otimes K_{j_l},
	\end{equation}
	where $\sum_{s = 1}^l (i_s + j_s) = m$ for some $m \in \mathbb{N}$, $i_1 \geq 0 $, and $ j_k \geq 1$, for all $ 1\leq k \leq l-1 $. By Proposition~\ref{prop:fibres_product}, the vector space in \eqref{eq:component} is a direct summand of $(\mathcal{H} \star \mathcal{K})_m$, which implies that $\mathcal{F}(\mathcal{H}) \ast \mathcal{F}(\mathcal{K}) \subset 
    \mathcal{F}(\mathcal{H} \star \mathcal{K})$. 
	
	On the other hand, the reverse inclusion follows from the fact that each summand in $(\mathcal{H} \star \mathcal{K})_m$ is of the form $H_{i_1} \otimes K_{j_1} \otimes H_{i_2} \otimes \dots \otimes H_{i_l} \otimes K_{j_l}$, with $\sum_{s = 1}^l (i_s + j_s)= m$, hence a summand in the free product $\mathcal{F}(\mathcal{H}) \ast \mathcal{F}(\mathcal{K})$.
\end{proof}

By induction and by associativity of the operations of free product of Hilbert spaces and subproduct systems, the claim holds for finitely many quadratic subproduct systems:
\begin{cor} \label{cor: free product sps}
Let $\mathcal{H}^{(i)}$, for $i=0,\dots.,n$, be $n+1$ quadratic subproduct systems of finite-dimensional Hilbert spaces. The Fock space $\mathcal{F}(\star_{i = 0}^n \mathcal{H}^{(i)})$ is unitarily isomorphic to the Hilbert space free product of the Fock spaces $\ast_{i = 0}^n \mathcal{F}({\mathcal{H}^{(i)}})$.   
\end{cor}

\section{Toeplitz algebras and KK-theory}
\label{sec:Toe}
\subsection{Free products of Toeplitz algebras and functoriality}
Note that the free product of algebras in \eqref{eq:free_prod_alg} is \emph{the} coproduct in the category of associative algebras over a field. For \Cs algebras, there are minimal and maximal free products, making the question of which free product is the right categorical coproduct particularly relevant. As discussed in \cite{Avi82}, if one considers GNS representations together with designated cyclic vectors, one can define a "free product representation" with a
designated cyclic vector, thus obtaining what Avitzour \cite{Avi82} calls a \emph{small free product representation}. Let us recall how the two constructions work and relate to each other.

\begin{definition}
    Given two separable and unital C*-algebras $A_1$ and $A_2$, their \emph{unital full free product} is given by the following commuting diagram of one-to-one unital morphisms:
\[ \xymatrix{
\mathbb{C}\ar[d] \ar[r] & A_1  \ar[d]
\\ A_2 \ar[r] & A_1 \star A_2 }.
\]
\end{definition}

\begin{definition} \label{defn:redFP}
	Let $\{(A_i, \phi_i, \mathcal{H}_{\phi_i}, \xi_i): i=0, \dots, n\}$ be a family of unital $C^*$-algebras with GNS states, Hilbert spaces, and unit vectors. Let $\lambda_i$ denote the left multiplication. The reduced free product $(A, \phi) = \ast^{r}_i (A_i, \phi_i)$ is the $C^*$-subalgebra generated by $\cup_{i \in I} \lambda_i(A_i)$, in the free product Hilbert space $\star_{i=0}^n (\mathcal{H}_{\phi_i}, \xi_i)$.
\end{definition}

Consider now the Toeplitz algebra $\Toe_{\mathcal{H}}$ of a subproduct system of Hilbert spaces $\mathcal{H}$. By construction, $\Toe_{\mathcal{H}}$ acts on the Fock space $\mathcal{F}(\mathcal{H})= \bigoplus_{m \geq 0} H_m$ faithfully via the left shift operators. We denote this $*$-representation by $(\mathcal{F}_{\mathcal{H}},\tau)$ and we refer to it as the Toeplitz representation (cf. \cite[Definition~2.13]{Vis11}).

\begin{prop}
\label{prop:state}
The Toeplitz representation $(\mathcal{F}_{\mathcal{H}},\tau)$ is equivalent to the GNS representation induced by the state $\varphi: \Toe_{\mathcal{H}} \to \C$ given by $\varphi(T) := \inp{T(\omega_0), \omega_0}$, with $\omega_{0}$ the unit vector in $H_0$.
\end{prop}
Note that the state $\varphi$ is the projection from $\Toe_{\mathcal{H}}$ onto the complex numbers. This fact, combined with the discussion in the previous section, yields:
\begin{theorem} \label{thm:toeplitz_free}
	Let $\mathcal{H}$ and $\mathcal{K}$ be quadratic subproduct systems, and let $\mathcal{H} \star \mathcal{K}$ be the free product of the subproduct systems in Definition \ref{defn:free_quad_SPS}. We have $$\Toe_{\mathcal{H}} \ast_\C \Toe_{\mathcal{K}} \cong \Toe_{\mathcal{H} \star \mathcal{K}}.$$
\end{theorem}

We shall now elaborate more on the categorical aspects of our construction.  Recall that the category $\Hilbf$ of finite-dimensional Hilbert spaces is the primary example of a strict \Cs tensor category, with morphisms being linear maps and the unit object being the uniquely defined one-dimensional Hilbert space $\mathbb{C}$. In particular, we consider the subcategory $\Hilf$ where morphisms are isometric maps.

Recall that a \emph{lax monoidal functor} \cite[Section~XII.2]{SML} is a functor between two monoidal categories,
together with two coherence maps satisfying associativity and unitality.
\begin{definition}
A subproduct system is a lax monoidal functor from $(\mathbb{N}_{0}, +,0)$ to $\Hilf$. 
The category of subproduct systems of finite-dimensional Hilbert spaces $\SPSf$ is the subcategory of $ (\Hilf)^{\mathbb{N}_{0}}$ whose objects are lax monoidal functors from $\mathbb{N}_{0}$ to $\Hilf$ and whose morphisms are as in \cite[Definition~ 1.4]{ShSo09}.
\end{definition}

There is a functor $\mathbf{Toe}$ from the category $\SPSf$ to the category of unital separable C*-algebras with states that associates to every subproduct system of Hilbert spaces $\mathcal{H}$ the corresponding Toeplitz C*-algebra with distinguished state $(\mathbb{T}_{\mathcal{H}},\tau)$ as in \ref{prop:state}, and to every morphism of subproduct systems the corresponding *-homomorphism at the level of \Cs algebras.

Our construction of the free product of subproduct system is therefore natural, as it is mapped to the corresponding free product of \Cs algebras by the functor $\mathbf{Toe}$.

\begin{example}

The Cuntz algebras $\mathcal{O}_n$ can be realised as quotients of the $n$-fold free product of Toeplitz algebras $\Toe$, which in turn can be realised as Toeplitz algebras of the (sub)product system with
$H_n:=\mathbb{C}$ for all $n$. 
    
\end{example}

\begin{example}[Cuntz--Krieger algebras and quadratic monomial ideals]
\label{ex:CK}
Monomial ideals are a special class of ideals, and subproduct systems associated with monomial ideals \cite{KakSha} give rise to many well-studied operator algebras, including Cuntz–-Krieger and subshift C*-algebras \'a la Matsumoto \cite{KM97}. 

Let us recall how Cuntz--Krieger algebras can be described using subproduct systems. Let $A \in \mathrm{Mat}_{n}\lbrace 0,1 \rbrace $, with no row or column equal to zero.  Inside the free algebra $\mathbb{C} \langle X_1, \dots, X_n \rangle$ we consider the quadratic monomial ideal generated by the quadratic monomials corresponding to the zero entries of the matrix $A$, i.e.,
\[ J_A := \langle X_iX_j \ : \ 1 \leq i, j \leq n, A_{ij}=0 \rangle . \]
The corresponding subproduct system $\mathcal{H}^{A}$ has fibres $H^{A}_{m}$ spanned by the admissible words of length $m$ in the alphabet $\lbrace X_1, \dots, X_n \rbrace$ and agrees with the subproduct system of a Matsumoto shift in the sense of \cite{ShSo09}. Note that $A$ can be interpreted as the incidence matrix of the underlying Markov chain. 

Let now $B \in  \mathrm{Mat}_{m}\lbrace 0,1 \rbrace $, satisfying again the condition of having no row or column equal to zero. Consider the associated subproduct system $\mathcal{H}^{B}$ and take the free product subproduct system $\mathcal{H}^A \star \mathcal{H}^B$. It is easy to see that this is the subproduct system of the Cuntz--Krieger algebra of the $(n+m) \times (n+m)$ matrix
\begin{equation*}
        \begin{pmatrix}
   A & \vline &  \mathbf{1}\\
\hline
  \mathbf{1} & \vline & B
\end{pmatrix},\end{equation*}
where $\mathbf{1}$ denotes the matrix with all entries equal to one. 

To describe what this means at the level of the underlying Markov chains, we need the following standard notion from graph theory. 
\begin{definition}[{\cite[Page~21]{FH69}}]
Given two directed graphs $E = (E^0, E^1), F = (F^0, F^1)$, their join $E + F$ is the graph with vertex set 
\[ (E + F)^0 := E^ 0 \cup E^1, \] and edge set \[ (E+F)^{1} := E^1 \cup F^1 \cup \left\{ (v_1, v_2), (v_2, v_1) : \forall v_1 \in E^0, v_2 \in F^0 \right\}, \] where $(v_1, v_2)$ denotes the arrow from $v_1$ to $v_2$. 
\end{definition}
Given matrices $A$ and $B$ with underlying Markov chains $E_A$ and $E_B$, the free product of their subproduct systems is the subproduct system with underlying Markov chain their graph join $E_A + E_B$. Correspondingly, the Cuntz--Pimsner algebra is the Cuntz--Krieger algebra of the incidence matrix of the graph join of the two underlying directed graphs.
\end{example}

 The examples above can also be interpreted as special cases of a result on Cuntz--Pimsner algebras of correspondences, due to Speicher \cite{Sp98} (see also \cite{BrOz}), stating that the Toeplitz algebra of the finite direct sum of C*-correspondences over the same coefficient algebra $A$ is the amalgamated free product over $A$ of the Cuntz--Pimsner algebras of a single correspondence.

 \begin{theorem}[{\cite[Example~4.7.5]{BrOz}}]
\label{thm:Sp}
  Let $H_{i}$, be a family of C*-correspondences over $A$. Denote the corresponding Toeplitz--Pimsner algebras by $\mathcal{T}(H_i)$, and by $E_{H_i}$ their conditional expectations. We have
\[ \left( \mathbb{T} ( \oplus_{i}H_{i} ), E_{\oplus H_i} \right) \simeq  \star^r_{i} (\mathbb{T}(H_i), E_{H_i}).\]
\end{theorem}

In the next section, we will focus on free products of a special class of \emph{proper} subproduct systems. To our knowledge, this is the first time that such examples have been studied. However, before doing that, we shall first discuss the KK-theory and nuclearity of our Toeplitz algebras.

\subsection{KK-theory and nuclearity for free products}
\label{ss:KKtheory}

We shall now recall some known results about the K-theory of free products of \Cs algebras. Our main references are \cite{Avi82,EG2,Ger:Crelle}. 

    For nuclear $C^*$-algebras $A_1$ and $A_2$ and any separable $C^*$-algebra $E$, Germain has proven that the reduced free product $A_1 \star_\C A_2$ is KK-equivalent to the unital full free product $A_1 \star A_2$. This implies the existence of the following six-term exact sequence:
    \begin{equation}
    \begin{tikzcd}
    \label{eq:6tes_red}
    KK_0(E, \mathbb{C}) \ar{r}{(i_1+i_2)_{\star}} & KK_0(E, A_1 \oplus A_2) \ar{r}{(j_1-j_2)_{\star}} & KK_0 \ar{d}(E, A_1\star_\C A_2) \\
    KK_1(E, A_1 \star_\C A_2) \ar{u} & KK_1 \ar{l}(E, A_1 \oplus A_2) \ar{l}{(j_1 - j_2)_{\star}} & KK_1(E, \mathbb{C}) \ar{l}{(i_1 + i_2)_{\star}}.
    \end{tikzcd} 
\end{equation}

%\todo[inline]{Something more needs to be checked.}

We will use \eqref{eq:6tes_red} to compute the $K$-theory groups of the Toeplitz algebra of the reduced free product of two quadratic subproduct systems.
%\end{remark}

\begin{definition}[{\cite[Definition~5.1]{EG2}}]
    A unital $C^*$-algebra $A$ is said to be \emph{$K$-pointed} if there exists $\alpha \in KK(A, \C)$ such that $i_A^*(\alpha) = 1_\C$ with $i_A$ the inclusion of $\C$ in $A$ given by the unit.
\end{definition}
Observe that if $A$ is a unital $C^*$-algebra and $KK$-equivalent to $\mathbb{C}$, then $A$ is $K$-pointed.

\begin{theorem}[{\cite[Theorem~5.5]{EG2}}]\label{thm:K-pointed}
    Let $A_0$ and $A_1$ be two $K$-pointed $C^*$-algebras. Then $A_0 \oplus A_1$ is $KK$-equivalent to $A_0 \star A_1 \oplus \C$. 
\end{theorem}

\begin{cor} \label{cor:KK-equi}
    Let  $A_0$  and  $A_1$  be two  $K$-pointed  $C^*$-algebras. If  $A_0$  and  $A_1$  belong to the UCT class  $\mathcal{N}$, then their reduced free product  $A_0 \star_\C A_1$  also belongs to $\mathcal{N}$.
\end{cor}
\begin{proof}
The UCT class $\mathcal{N}$ is closed under direct sum and $KK$-equivalence. This fact, combined with Theorem \ref{thm:K-pointed}, yields $A_0 \oplus A_1 \sim_{KK} A_0 \star A_1 \oplus \C \in \mathcal{N}$. Furthermore, since  $A_0 \star A_1$  is  $K$-dominated by  $A_0 \star A_1 \oplus \mathbb{C}$, the claim follows. 
    %Theorem \ref{thm:K-pointed} implies that if $A_0$ and $A_1$ belong to the UCT class $\mathcal{N}$, then $(A_0 \star A_1 )\oplus \C$ also belongs to the UCT class $\mathcal{N}$. Since $A_0 \star A_1$ is $K$-dominated by $(A_0 \star A_1 )\oplus \C$, we have $A_0 \star A_1$ belongs to the UCT class $\mathcal{N}$. Therefore, to establish that $A_0 \star A_1$ is $KK$-equivalent to $\C$, it suffices to show that $A_0 \star A_1$ has the same $K$-theory groups as $\C$ since the UCT ensures that $KK$-equivalence is determined by $K$-theory when $C^*$-algebras belong to $\mathcal{N}$.
\end{proof}

\begin{theorem} \label{thm:fp_kk-eq}
    Let $\mathcal{H}$ and $\mathcal{K}$ be standard subproduct systems of Hilbert spaces. Assume that the Toeplitz algebras $\Toe_{\mathcal{H}}$ and $\Toe_{\mathcal{K}}$ are both \emph{nuclear} and KK-equivalent to the complex numbers. Then so is the Toeplitz algebra $\Toe_{\mathcal{H}\star \mathcal{K}}$.
\end{theorem}
\begin{proof}
By Corollary \ref{cor:KK-equi}, under the assumptions of the theorem, it follows that $\Toe_H \star_\C \Toe_K \in \mathcal{N}$. Consequently, it suffices to show that $\Toe_{\mathcal{H}} \star \Toe_{\mathcal{K}}$ has the same $K$-theory as $\mathbb{C}$, thanks to the fact that $\Toe_{\mathcal{H} \star \mathcal{K}} \cong \Toe_{\mathcal{H}} \star_\C \Toe_{\mathcal{K}}$. %By assumption, $\Toe_{\mathcal{H}}, \Toe_{\mathcal{K}}$ are nuclear and $KK$-equivalent to $\C$. 
After replacing  $E$ with $\mathbb{C}$ and  $A_0, A_1$ with $\Toe_{\mathcal{H}}$ and  $\Toe_{\mathcal{K}}$, respectively, the long exact sequence \eqref{eq:6tes_red} becomes
\[
0 \longrightarrow K_1(\Toe_{\mathcal{H}} \star_\C \Toe_{\mathcal{K}}) \longrightarrow \mathbb{Z} \xrightarrow{(i_1 + i_2)_\ast} \mathbb{Z} \oplus \mathbb{Z} \longrightarrow K_0(\Toe_{\mathcal{H}} \star_\C \Toe_{\mathcal{K}}) \longrightarrow 0.
\]
From this, we compute that 
\begin{align*}
    &K_1(\Toe_{\mathcal{H}} \star_\C \Toe_{\mathcal{K}}) \cong \ker{(i_1 + i_2)_\ast} \cong \{ 0 \} \cong K_1(\C), \\
    &K_0(\Toe_{\mathcal{H}} \star_\C \Toe_{\mathcal{K}}) \cong \text{coker}{(i_1 + i_2)_\ast} \cong \mathbb{Z} \cong K_0(\C),
\end{align*}
which implies that $\Toe_{\mathcal{H}} \star_\C \Toe_{\mathcal{K}}$ is $KK$-equivalent to $\C$.
By the definition of free product of subproduct systems and functoriality, we have $\Toe_{\mathcal{H}} \star_\C \Toe_{\mathcal{K}}\cong \Toe_{\mathcal{H} \star \mathcal{K}}$. Therefore, $\Toe_{\mathcal{H} \star \mathcal{K}}$ is $KK$-equivalent to $\C$.

By assumption, $\Toe_{\mathcal{H}}, \Toe_{\mathcal{K}}$ are nuclear and $KK$-equivalent to $\C$. Moreover, the compact operators $\mathbb{K}(\mathcal{F})$ are contained in the Toeplitz algebra $\Toe_{\mathcal{H}}$ (see. \cite[Corollary~3.2]{Vis12}). Thus, thanks to \cite[Theorem~1.1]{Oz02}, we obtain nuclearity of the reduced free product $\Toe_{\mathcal{H}} \star_\C \Toe_{\mathcal{K}}$.
\end{proof}

\section{A case study: free products of Temperley--Lieb subproduct systems} \label{sec:TL}
\begin{definition}[{\cite[Definition~1.2]{HaNe21}}]\label{def:TP_pol}
Let $H$ be a finite-dimensional Hilbert space of dimension $m \geq 2$. A non-zero vector $P \in H\otimes H$ is called \emph{Temperley--Lieb} if there is $\lambda > 0$ such that the orthogonal projection $e\colon H\otimes H \to \C \cdot P$ satisfies
\[ (e\otimes 1)(1\otimes e)(e\otimes 1) = \dfrac{1}{\lambda}(e\otimes 1)\quad\text{in}\quad B(H\otimes H\otimes H).\]
\end{definition}

The standard subproduct system $\mathcal{H}^P$ defined by the ideal $\langle P\rangle\subset T(H)$ generated by $P$ is called a Temperley--Lieb subproduct system. We write $\mathcal{F}_P = \mathcal{F}_{\mathcal{H}_P}$, $\Toe_P = \Toe_{\mathcal{H}_P}$ and $\Pim_P = \Pim_{\mathcal{H}_P}$.

We will often fix an orthonormal basis in $H$ and identify $H^{\otimes n}$ with the space of homogeneous noncommutative polynomials of degree $n$ in variables $X_1,\dots,X_d$. In particular, we write a vector $P\in H\otimes H$ as a noncommutative polynomial $P=\sum^d_{i,j=1}a_{ij}X_iX_j$. Consider the matrix $A=(a_{ij})_{i,j}$. By \cite[Lemma~1.4]{HaNe21}, $P$ is Temperley--Lieb if and only if the matrix $A\bar A$ is unitary up to a (non-zero) scalar factor, where $\bar A=(\bar a_{ij})_{i,j}$. Since the ideal generated by $P$ does not change if we multiply $P$ by a non-zero factor, we may always assume that $A\bar A$ is unitary.

The following result gives a complete set of relations in $\Toe_P$.

\begin{theorem}[{\cite[Theorem~2.11]{HaNe22}}] \label{thm: univ_relations}
Let $A=(a_{ij})_{i,j}\in \mathrm{GL}_d(\C)$ ($m\geq 2$) be such that $A\bar A$ is unitary. Let $q\in(0,1]$ be the number such that $\mathrm{Tr}(A^*A)=q+q^{-1}$. Consider the noncommutative polynomial $P=\sum^d_{i,j=1}a_{ij}X_iX_j$.  Then $\Toe_P$ is a universal C$^*$-algebra generated by the \Cs-algebra $c := C(\mathbb{Z}_+ \cup \{\infty\})$ and elements $S_1,S_2,...,S_d$ satisfying the relations
\[ fS_i = S_i\gamma(f) \quad (f \in c,\ 1 \leq i \leq d), \quad \sum^d_{i=1} S_iS_i^* = 1 - e_0, \quad \sum^d_{i,j=1} a_{ij}S_iS_j = 0,  \]
\[S_i^*S_j + \phi\sum^d_{k,l=1}a_{ik}\bar{a}_{jl}S_kS_l^* = \delta_{ij}1 \quad (1 \leq i,j \leq d), \]
where $\gamma\colon c\to c$ is the shift to the left (so $\gamma(f)(n)=f(n+1)$), $e_0$ is the characteristic function of $\{0\}$ and $\phi\in c$ is the element given by
\begin{equation} \label{eq:phi(n)}
\phi(n)= \dfrac{[n]_q}{[n+1]_q},\quad\text{with}\quad [n]_q= \dfrac{q^{n}-q^{-n}}{q - q^{-1}}.
\end{equation}
\end{theorem}

Here $c$ is identified with a unital subalgebra of $\mathbb{K}(\mathcal{F}_P)+\C1 \subset \Toe_P$, with $e_n\in c$ being identified with the projection $\mathcal{F}_P \to H_n$. Note also that $\phi(n)\to q$ as $n\to+\infty$.

The relations become slightly simpler if we write $P$ in a standard form. Namely, by \cite[Proposition~1.5]{HaNe21}, up to a unitary change of variables and rescaling, we may assume that our Temperley--Lieb polynomial $P$ has the form 
\begin{equation}
    \label{eq:TLPoly_normal} P = \sum_{i=1}^{m} a_i X_iX_{m-i+1},    \quad \text{with}\quad \vert a_i a_{m-i+1} \vert =1.
\end{equation} 
%Note that this yields inner automorphisms $\alpha_{P,P'}:\Toe_P\to \Toe_{P'}$ and $\beta_{P,P'}:\Pim_P\to \Pim_{P'}$.

\subsection{Fusion rules for free products of Temperley--Lieb subproduct systems}
For $i=1,2$, let $P_i$ be Temperley--Lieb polynomials with associated Temperley--Lieb subproduct systems $\mathcal{H}^{P_i}$. For ease of notation, denote their free product $\mathcal{H}^{\langle P_1 , P_2 \rangle}$ by $\mathcal{H}^{\star}$. In particular, we have 
\begin{align*}
H_0^{\star} & \cong \mathbb{C};\\
H_1^{\star} & \cong H_1^{P_1} \oplus H_1^{P_2} \cong \mathrm{Span}_\C \{ e_1^1, \dots, e_{d_1}^1, e_1^2, \dots, e_{d_2}^2 \};\\
H_2^{\star} & \cong H_1^{\star} \otimes H_1^{\star} \ominus \mathrm{Span}_\C \left\{ P_1(e_1^1, \dots, e_{d_1}^1), P_2(e_1^2, \dots, e_{d_2}^2) \right\},
\end{align*}
where $d_i = \dim(H_1^{P_i})$, and $\{ e_j^i : j = 1, 2, \dots, d_i \}$ forms an orthonormal basis of $H_1^{P_i}$.

In \cite{HaNe22}, the author constructed maps
\[ w_n^{P} := ([2]_q \phi(n+1))^{\frac{1}{2}} (f_{n+1} \otimes 1)(1^{\otimes n}\otimes v) : H_{n}^{P} \to H_{n+1}^{P}\otimes H_1^{P}, \]
where $v = P(e) \in H_1^{P} \otimes H_1^{P}$ is the Temperley--Lieb vector corresponding to the polynomial $P$.

We should use a small adaptation of those maps to construct an explicit isometry $$w_n^{\star} : ({H}^{\star}_{n-1})^{\oplus 2} \oplus {H}^{\star}_{n+1} \to {H}^{\star}_{n} \otimes  {H}^{\star}_1.$$

For $i=1,2$, consider the two projections $Q_i: H_1^{\star} \cong H_1^{p_1} \oplus H_1^{p_2}$ onto $H_1^{p_1}$ and $H_1^{p_2}$, respectively.
\begin{equation}
\begin{split}
	 & w_n^i = ([2]_q \phi^i(n+1))^{\frac{1}{2}} Q_i^{\otimes n + 2} (f_{n+1} \otimes 1)(1^{\otimes n} \otimes v_i) Q_i^{\otimes n} \\
	& \, + ([2]_q\phi^i(n))^{\frac{1}{2}} (1 - Q_i) \otimes Q_i^{\otimes n+1} (1 \otimes f_n \otimes 1)(1\otimes 1^{\otimes n-1}\otimes v_i) (1 - Q_i) \otimes Q_i^{\otimes n-1} \\
	& \,+ ([2]_q\phi^i(n-1))^{\frac{1}{2}} 1\otimes (1 - Q_i) \otimes Q_i^{\otimes n} (1^{\otimes 2} \otimes f_{n-1} \otimes 1)(1^{\otimes n}\otimes v_i) 1\otimes(1 - Q_i) \otimes Q_i^{\otimes n-2} \\
%	&+ ([2]_q\phi^i(n-2))^{\frac{1}{2}} 1^{\otimes 2}\otimes (1 - Q_i) \otimes Q_i^{\otimes n-3} (1^{\otimes 3} \otimes f_{n-2} \otimes 1)(1^{\otimes 3} \otimes 1^{\otimes n-3}\otimes v_i) 1^{\otimes 2}\otimes(1 - Q_i) \otimes Q_i^{\otimes n-3} \\
%	&+ \quad \vdots \\
	& \, + \cdots + ([2]_q\phi^i(1))^{\frac{1}{2}} 1^{\otimes n-1}\otimes (1 - Q_i) \otimes Q_i^{\otimes 2} (1^{\otimes n} \otimes f_{1} \otimes 1)(1^{\otimes n}\otimes v_i) 1^{\otimes n-1}\otimes (1 - Q_i).  \label{decomp_free}
\end{split}
\end{equation}

\begin{prop}
	The map $w_n^{\star} := w_n^1 \oplus w_n^2 : ({H}^{\star}_{n-1})^{\oplus 2} \to{H}^{\star}_{n} \otimes  {H}^{\star}_1 $ is an isometry.
\end{prop}

Combining this with the structure map $\iota^{\star}_{n,1}: {H}^{\star}_{n+1} \mapsto {H}^{\star}_{n} \otimes {H}^{\star}_{1} $ we obtain the following:
\begin{theorem}
\label{thm:WRn}
There is a unitary isomorphism:
\[ W^R_n: = (w_n^{\star}, \iota_{n,1}^{\star}) : ({H}^{\star}_{n-1})^{\oplus 2} \oplus {H}^{\star}_{n+1}\to {H}^{\star}_{n} \otimes  {H}^{\star}_1 \]
\end{theorem}
\begin{proof}
	Let $v_i = P_i(e) \in H_1^{(i)} \otimes H_1^{(i)}$ be the Temperley--Lieb vector corresponding to the polynomial $P_i$. Using~\eqref{eq:freeprod}, it is not hard to see that $w_n^\star$ is orthogonal to $\iota_{n, 1}^\star$.
    
    Note that, in the image of the map $w_n^i$, the last component belongs to $H_1^{(i)}$, one of the orthogonal components. This yields \[\inp{w_n^i(\xi_1), w_n^j(\xi_2)} = 0, \quad \mbox{for all } \xi_1, \xi_2 \in H^{\star}_{n-1}, i \neq j.\] 

	It remains to prove that each $w_n^i$ is an isometry. Let $\xi \in H^\star_n$. It is not hard to see that the summands of $w_n^i(\xi)$ are mutually orthogonal, since the images of $Q_i$ and $1-Q_i$ are. 
    
    For any vector $\xi \in \text{Im}(1^{\otimes k}\otimes(1 - Q_i) \otimes Q_i^{\otimes n-k+1})$, we have that $\inp{w_n^k(\xi), w_n^k(\xi)}$ is equal to
	\begin{align*}
		\left\|([2]_q\phi^i(n-k))^{\frac{1}{2}} (1^{\otimes k} \otimes f_{n-k+1} \otimes 1)(1^{\otimes k} \otimes 1^{\otimes n-k}\otimes v_i)(\xi) \right\|^2.
	\end{align*}
	Write $(1^{\otimes k} \otimes f_{n-k+1} \otimes 1)(1^{\otimes k} \otimes 1^{\otimes n-k}\otimes v_i)$ as $1^{\otimes k} \otimes (f_{n-k+1}\otimes 1)(1^{\otimes n-k}\otimes v_i)$, the tensor product of linear operators. Let $\{ \eta_j : j = 1, 2, \dots, \dim(H^\star_1)^k \}$ be the basis of ${H^\star_1}^{\otimes k}$ and write $\xi = \sum_{j = 1}^{\dim(H_1^\star)^k} \eta_j \otimes \zeta_j$ for some $\zeta_j \in H^\star_{n-k}$, then we have 
	\begin{align*}
		&\left\|([2]_q\phi^i(n-k))^{\frac{1}{2}} (1^{\otimes k} \otimes f_{n-k+1} \otimes 1)(1^{\otimes k} \otimes 1^{\otimes n-k}\otimes v_i)(\xi) \right\|^2 \\
		&=  \left\|([2]_q\phi^i(n-k))^{\frac{1}{2}} 1^{\otimes k} \otimes (f_{n-k+1}\otimes 1)(1^{\otimes n-k}\otimes v_i)\left(\sum_{j = 1}^{\dim(H^\star_1)^k} \eta_i \otimes \zeta_i\right) \right\|^2 \\
		&= \left\|\sum_{j=1}^{\dim(H^\star_1)^k} \eta_i \otimes ([2]_q\phi^i(n-k))^{\frac{1}{2}} \cdot (f_{n-k+1}\otimes 1)(1^{\otimes n-k}\otimes v_i) (\zeta_i) \right\|^2 \\
		&= \sum_{j = 1}^{\dim(H^\star_1)^k} \|\zeta_j\|^2 = \| \xi \|^2,
	\end{align*}
	where the second to last equality follows from {\cite[Equation~4.1]{HaNe22}}. This proves the claim.
\end{proof}

For ease of notation, we write 
\begin{align*}
    &V_R := \oplus_{n = 1}^\infty W_n^R : \oplus_{n = 0}^\infty (H_n^\star)^{\oplus 2} \cong \mathcal{F}^{\oplus 2} \to \oplus_{n = 1}^\infty H_n^\star \otimes H_{1}^\star \cong \mathcal{F}_+ \otimes H_{1}^\star, \\
   &w^i = \oplus_{n = 0}^\infty w_n^i : \oplus_{n = 0}^\infty H_n^\star \cong \mathcal{F} \to \oplus_{n = 1}^\infty H_n^\star \otimes H_{1}^\star \cong \mathcal{F}_+ \otimes H_{1}^\star, \text{ for each $i = 1, 2$.}
\end{align*}
These two maps are important in the construction of $KK$-element in $KK(\Toe, \C)$ in Section~\ref{section: An explicit KK-equivalence}.
\begin{example}
	Consider the Temperley--Lieb polynomials $$P_1(X_1, X_2) = X_1X_2 - X_2X_1,\quad P_2(X_3, X_4, X_5) = a_1 \cdot X_3X_5 + a_2 \cdot X_4X_4 + a_3 \cdot X_5X_3, $$ such that $|a_1|^2 + |a_2|^2 + |a_3|^2 = q + q^{-1}$ for $q \in (0, 1]$.
 The Temperley--Lieb subproduct systems induced by $P_1 $ and $P_2$ are the maximal subproduct systems with fibres
 \[\begin{aligned}
     H_1^{P_1} & = \mathrm{Span}_\C \{ e_1, e_2 \}, \\
     H_2^{P_1}  & = \mathrm{Span}_\C \lbrace e_1\otimes e_2 - e_2 \otimes e_1\rbrace^\perp,\end{aligned}\]
     and
\[\begin{aligned} H_1^{P_2} & = \mathrm{Span}_\C\{ e_3, e_4, e_5 \}, \\   H_2^{P_2}  & = \mathrm{Span}_\C \lbrace a_1 \cdot e_3\otimes e_5 + a_2 \cdot e_4 \otimes e_4 + a_3 \cdot e_5 \otimes e_3\rbrace^\perp, \end{aligned}\]
respectively. 
 
The free product subproduct system has fibres 
\[\begin{aligned}
    H^{\langle P_1,P_2\rangle}_1 & = H_1^{P_1} \oplus H_1^{P_2}, \\
H^{\langle P_1,P_2\rangle}_2 & = \mathrm{Span}_\C\{ e_1\otimes e_2 - e_2 \otimes e_1, a_1 \cdot e_3\otimes e_5 + a_2 \cdot e_4 \otimes e_4 + a_3 \cdot e_5 \otimes e_3\}^\perp, 
\end{aligned}\]
and satisfies \eqref{eq:decomp_free_prod}. Consider the case when $n = 1$, we have the decomposition:
 \[(w_1^1, w_1^2) :H^{\langle P_1,P_2\rangle}_1 \oplus H^{\langle P_1,P_2\rangle}_1 \to H^{\langle P_1,P_2\rangle}_2 \otimes H^{\langle P_1,P_2\rangle}_1\] 
 
Applying $w_1^2$ defined in \eqref{decomp_free}, we obtain:
\begin{align*}
    &w_1^2(e_1) = e_1 \otimes v_2, \\
    &w_1^2(e_2) = e_2 \otimes v_2, \\
    &w_1^2(e_3) =  \frac{[2]_q}{[3]_q^{1/2}}  \cdot \left(e_3 \otimes v_2 - \frac{a_1a_3}{q + q^{-1}} \cdot v_2 \otimes e_3 \right),\\
    &w_1^2(e_4) = \left( \frac{[2]_q^2}{[3]_q} \right)^{\frac{1}{2}}\cdot \left( e_4 \otimes v_2 - \frac{a_2^2}{q + q^{-1}}\cdot v_2 \otimes e_4 \right),\\
    &w_1^2(e_5) = \left( \frac{[2]_q^2}{[3]_q}\right)^{\frac{1}{2}} \cdot \left(e_5 \otimes v_2 - \frac{a_1a_3}{q + q^{-1}}\cdot v_2 \otimes e_5 \right),
\end{align*}
where $v_2 := a_1 \cdot e_3\otimes e_5 + a_2 \cdot e_4 \otimes e_4 + a_3 \cdot e_5 \otimes e_3$ is the quadratic relation that defines Temperley--Lieb subproduct system $\mathcal{H}^{P_2}$.

\end{example}
By associativity, the above construction can be extended to the free products of a finite number of Temperley--Lieb subproduct systems:
\begin{theorem}
    Let $P_i, i = 1, \dots, r$ be $r$ Temperley--Lieb polynomials with associated Temperley--Lieb subproduct systems $\mathcal{H}^{P_i}$. Denote by $\mathcal{H}^\star$ their free product. There is a unitary isomorphism 
    \[ 
        W^R_n := (w^\star, \iota_{n, 1}^\star) : (H^\star_{n-1})^{\oplus r} \oplus {H}^\star_{n + 1} \to {H}^\star_n \otimes {H}_1^\star,
    \]
    where $\iota_{n, 1}^\star$ are the structure maps of $\mathcal{H}^\star$, and $w^\star := (w_n^1, w_n^2, \dots, w_n^r)$ with $w_n^i$ as in \eqref{decomp_free}.
\end{theorem}

\subsection{Gysin Sequences}
\label{sec:Gysin}
In this section, we will construct the noncommutative Gysin sequence for the Toeplitz algebra of the free product of Temperley--Lieb subproduct systems. This will allow us to simplify the six-term exact sequence in $K$-theory induced by the extension \eqref{eq:CPext}. 

We start by recalling what is known about the $K$-theory of the Toeplitz and Cuntz--Pimsner algebras of a subproduct system.

\begin{theorem}[{\cite[Theorem~3.1 and Corollary~4.4]{HaNe22}}]
For a Temperley–Lieb polynomial~$P$, the inclusion $i: \mathbb{C} \to \Toe_P$ is a $KK$-equivalence. Moreover, we have $[e_0]=(2-m)[1]$ in~$K_0(\Toe_P)$.
\end{theorem}

As a consequence, the six-term exact sequence induced by the defining extension \eqref{eq:CPext} simplifies notably, and one obtains the following result about the $K$-theory of the Cuntz--Pimsner algebra of a Temperley--Lieb subproduct system. 

\begin{cor}[{\cite[Corollary~4.4]{HaNe22}}]\label{cor:Ktheory}
For every Temperley--Lieb polynomial $P$ in $d$ variables, 
    \[K_0(\Pim_{P}) \cong \mathbb{Z}/(d-2)\mathbb{Z}, \qquad K_1(\Pim_{P}) \cong \begin{cases}
        \mathbb{Z}, & d=2,\\
        0, & d \geq 3.
    \end{cases}\]
\end{cor} 

Given that the Toeplitz algebras associated with Temperley--Lieb subproduct systems satisfy the conditions of Theorem~\ref{thm:fp_kk-eq}, we derive the following result:
\begin{theorem}
    Let $P_1$ and $P_2$ be Temperley--Lieb polynomials, with Toeplitz algebras $\Toe_{P_1}$ and $\Toe_{P_2}$. Then the Toeplitz algebra of the free product subproduct system $\Toe_{\langle P_1, P_2 \rangle}$ is isomorphic to the reduced free product $\Toe_{P_1} \star \Toe_{P_2}$ and it is $KK$-equivalent to the algebra of complex numbers $\C$.
\end{theorem}
\subsubsection{An explicit KK-equivalence}
\label{section: An explicit KK-equivalence}
We shall now make our $KK$-equivalence result more explicit. %by defining a quasi-homomorphism $(\psi_-, \psi_+)$ which represents a $KK$-element in $KK(\Toe_{\langle P, Q \rangle}, \C)$. 
To do so, we shall employ arguments similar to those in \cite{AK21,HaNe21}.

With the same notation as in the previous section, let $P_1$ and $P_2$ be two Temperley--Lieb polynomials with associated subproduct systems $\mathcal{H}^{P_1}$ and $\mathcal{H}^{P_2}$, and let us consider their free product by $\mathcal{H}^{P_1} \star \mathcal{H}^{P_2} = \mathcal{H}^{\langle P_1, P_2 \rangle}.$
In what follows, we will omit the subscript $\langle P_1, P_2 \rangle$ and write $\Toe$ for $\Toe_{\langle P_1, P_2 \rangle}$ and $\mathcal{F}$ for $\mathcal{F}_{{\langle P_1, P_2 \rangle}}$. 

By Theorem~\ref{thm:WRn} we have maps $W^R_n: H_n \otimes H_1 \to H_{n+1} \oplus H_{n-1}^{\oplus 2}$, for every $n$. This allows us to construct a map $W_R := (\iota, V_R)^* : \mathcal{F} \otimes H_1 \to  \mathcal{F}^{\oplus 3}$, with range $\mathcal{F}_+ \oplus \mathcal{F}^{\oplus 2}$. 

%for every $n$ we have isometries $w_n : E_n^{\oplus r} \to E_{n+1} \otimes E_1$. Define the direct sum $V_R := \sum_{n = 0}^{\infty} w_n : \oplus_{n = 0}^\infty E_n^{\oplus r} \to \oplus_{n = 0}^\infty E_{n + 1} \otimes E_1 =: F_+ \otimes E_1$, where $F_+ = \oplus_{n = 1}^\infty E_n$. Denote the inclusion $\iota_n: E_n \hookrightarrow E_{n-1}\otimes E_1$, we have $w_{n-1}w_{n-1}^* + v_{n+1}v_{n+1}^* = 1$ due to fusion rules. Therefore, if we write $\iota := \sum_{n=1}^\infty \iota_n : \oplus_{n = 1}^\infty E_n \to E_{n-1} \otimes E_1$, we may define the adjoint of the decomposition  we have a map $W_R := (\iota, V_R)^* : F \otimes E_1 \to F \oplus F^{\oplus r}$, which generalizes {\cite[Equation~5.1]{AK21}} and {\cite[Lemma~4.1]{HaNe22}}, and admits range $F_+ \oplus F^{\oplus r}$. Note that this map will play an essential role in proving $KK$-equivalence in Section \ref{sec:Gysin}.

We consider the pair of homomorphisms $(\psi_+, \psi_-)$, where $\psi_\pm : \Toe \to \mathcal{L}(\mathcal{F}^{\oplus 3})$
\begin{align*}
    &\psi_+(x) = x^{\oplus 3}, \\
    &\psi_-(x) = W_R(x \otimes 1_{H_1})W_R^*.
\end{align*}
We will show that the above pair $(\psi_+, \psi_-)$ gives a $KK$-element  which is a left and right inverse to the $KK$-class of the inclusion $i: \C \to \Toe$. 
\begin{lemma}
	The pair $(\psi_+, \psi_-)$ defines an element $[\psi_-, \psi_+]$ in $KK(\Toe, \C)$.
\end{lemma}

\begin{proof}
	It is sufficient to prove that for all $x \in \Toe$, $\psi_+(x) - \psi_-(x) \in \mathbb{K}(\mathcal{F}^{\oplus 3})$. 
    
    Let $d_1 = \dim(H_1)$. Since the Toeplitz algebra $\Toe$ is generated by the Toeplitz operators $T_i^*:=T_{e_i}^*$, where $\{ e_i : i = 1, 2, \dots, d_1 \}$ is an orthonormal basis of $H_1$, we only need to show that $\psi_+(T_i^*) - \psi_-(T^*_i)$ is a compact operator for all $i$.

Writing $W_R(x \otimes 1_{H_1})W_R^*$ in matrix form, we are left with checking that	
\begin{align*}
	\begin{bmatrix}
		T_i^* - \iota^*(T_i^* \otimes 1) \iota & \iota^*(T_i^* \otimes 1) w^1 & \iota^*(T_i^* \otimes 1)w^2 \\
		(w^1)^*(T_i^* \otimes 1) \iota & T_i^* - (w^1)^*(T_i^* \otimes 1)w^1 & (w^1)^* (T_i^* \otimes 1)w^2 \\
		(w^2)^*(T_i^* \otimes 1) \iota & (w^2)^* (T_i^* \otimes 1)w^1 & T_i^* - (w^2)^*(T_i^* \otimes 1)w^2
	\end{bmatrix} = 0 \mod \mathbb{K}(\mathcal{F}^{\oplus 3}).
	\end{align*}
	Since $\iota$ is the inclusion of $H_{n+1}$ into $H_n \otimes H_1$ and $T_i^*$ commutes with the structure maps of the subproduct system, we have
	\[	(T_i^*\otimes 1)\iota(\xi) = T_i^*(\xi_k) \otimes e_k = \iota(T_i^* \otimes 1)	
\]
for all $\xi \in H_{n+1} \subset H_n \otimes H_1$, where we write $\xi = \sum_{k} \xi_k \otimes e_k$.  Consequently, \[ (T_i^* \otimes 1)\iota = \iota(T_i^* \otimes 1).\]
Moreover, $(w^i)^* \iota = 0$, so the lower triangular part of $\psi_+(T_i^*) - \psi_-(T_i^*)$ vanishes.
    
To finalize the proof, we need to show that $(T_i^* \otimes 1) w^k = w^kT_i^* \mod \mathbb{K}$, which together with the fact that $w^k$ is an isometry for all $k = 1, 2$, gives the desired result that $(\psi_+ - \psi_-)(T_i^*) = 0 $ modulo compact operators.
	
	Let $\xi \in H_n$, using the fact that $((T_i^* \otimes 1) w^k - w^kT_i^*)(\xi) = T^*_i(w_n - 1 \otimes w_{n-1})(\xi)$ as in the proof of {\cite[Lemma~4.1]{HaNe22}}, it suffices to prove that 
	\[ \lim_{n \to \infty} \|(w_n^k - 1 \otimes w_{n-1}^k) f_n \| = 0, \quad \mbox{for all \,} k = 1, 2. \]
	To this end, we have 
	\begin{align*}
		&((w_n^k - 1 \otimes w_{n-1}^k) f_n)^*((w_n^k - 1 \otimes w_{n-1}^k  ) f_n) \\
		&= f_n ((w_n^k)^*w_n^k - (w_n^k)^*(1 \otimes w_{n-1}^k) - (1\otimes w_{n-1}^k)^*w_n^k + 1 \otimes (w_{n-1}^k)^*w_{n-1}^k) f_n \\
		&= f_n(2 - (w_n^k)^*(1\otimes w_{n-1}^k) - (1\otimes w_{n-1}^k)^*w_n^k)f_n.
	\end{align*}
	Using the decomposition of $w_n^k$ from \eqref{decomp_free}, we compute,
	\begin{align*}
		(w_n^k)^*(1\otimes w_{n-1}^k) = \left( \frac{\phi(n)}{\phi(n+1)}\right)^{\frac{1}{2}}\cdot Q_k^{\otimes n} + \sum_{l = 0}^{n-1} 1^{\otimes l}\otimes (1 - Q_k) \otimes Q_k^{\otimes n-1-l},
	\end{align*}
	which can be written as a matrix
	\begin{align*}
		(w_n^k)^*(1\otimes w_{n-1}^k) &= \begin{bmatrix}
			\left( \frac{\phi(n)}{\phi(n+1)}\right)^{\frac{1}{2}} & 0 & 0 & \dots & 0 \\
			0 & 1 & 0 & \dots & 0 \\
			0 & 0 & 1 & \dots & 0 \\
			\vdots & \vdots & \vdots & \ddots & 0 \\
			0 & 0 & 0 & \dots & 0
		\end{bmatrix}, 
	\end{align*}
	acting on $\text{Im}(Q_k^{\otimes n}) \oplus \bigoplus_{l = 0}^{n-1}\text{Im}(1^{\otimes l}\otimes (1 - Q_k) \otimes Q_k^{\otimes n-1-l})$. 
    
    Applying the same reasoning to $(1\otimes w_{n-1}^k)^* w_n^k$ yields
	\begin{align*}
		f_n(2 - (w_n^k)^*(1\otimes w_{n-1}^k) - (1\otimes w_{n-1}^k)^*w_n^k)f_n = f_n\begin{bmatrix}
			2 - 2 \sqrt{\frac{\phi(n)}{\phi(n+1)}} & 0 & 0 & \dots & 0 \\
			0 & 0 & 0 & \dots & 0 \\
			0 & 0 & 0 & \dots & 0 \\
			\vdots & \vdots & \vdots & \ddots & 0 \\
			0 & 0 & 0 & \dots & 0
		\end{bmatrix}f_n,
	\end{align*}
	acting on $\mathrm{Im}(Q_k^{\otimes n}) \oplus \bigoplus_{l = 0}^{n-1}\mathrm{Im}(1^{\otimes l}\otimes (1 - Q_k) \otimes Q_k^{\otimes n-1-l})$.
    
Therefore, we obtain
	\begin{align*}
		\| f_n(2 - (w_n^k)^*(1\otimes w_{n-1}^k) - (1\otimes w_{n-1}^k)^*w_n^k)f_n \|^2 &= 2\left(1-\sqrt{\frac{\phi(n)}{\phi(n+1)}}\right)\\
		&= 2(1 - (1 - [n+1]_q^{-2})^{\frac{1}{2}}),
	\end{align*}
	which converges to zero as desired.
\end{proof}
In particular, we obtain the following \emph{explicit} KK-equivalence result.
\begin{theorem}\label{KK-equi-formula}
Let $P_1, P_2$ be two Temperley--Lieb polynomials, and let $\Toe_{\langle P_1, P_2 \rangle}$ be the associated free-product Toeplitz algebra. Denote by $i: \C \to \Toe_{\langle P_1, P_2 \rangle}$ the natural inclusion.
    The interior Kasparov product $[i]\otimes_{\Toe}[\psi_+, \psi_-]$ agrees with the unit $\textbf{1}_\C \in KK(\C, \C)$. In particular, $[i]$ and $[\psi_+, \psi_-]$ implement the $KK$-equivalence between $\C$ and $\Toe$. %$\Toe(L_{n_1} \oplus \dots \oplus L_{n_R})$. 
\end{theorem}
\begin{proof}
    The interior Kasparov product $[i]\otimes_{\Toe}[\psi_+, \psi_-] \in KK(\C, \C)$ is represented by the pair $(\psi_+\circ i, \psi_-\circ i)$, where $\psi_\pm \circ i: \C \to L(\mathcal{F}^{\oplus 3})$ are $*$-homomorphisms. 
    
    In particular, $\psi_+ \circ i$ is unital and \[\psi_- \circ i(1) = W_RW_R^*: \mathcal{F}^{\oplus 3} \rightarrow \mathcal{F}^{\oplus 3} \]is the orthogonal projection with range $\mathcal{F}_+ \oplus \mathcal{F}^{\oplus 2}$. 
    
    Therefore, $(\psi_+\circ i - \psi_-\circ i)(1) $ is the rank one orthogonal projection onto $\C \oplus \{0\} \subset \mathcal{F}^{\oplus 3}$. 

    Since we have already proven an abstract $KK$-equivalence between $\C$ and $\Toe$ in Theorem~\ref{thm:fp_kk-eq}, and $[i]$ maps the generators to generators, the claim follows.
\end{proof}
The defining extension of Cuntz--Pimsner algebras of subproduct systems of finite-dimensional Hilbert spaces \eqref{eq:CPext} induces a six-term exact sequence in $K$-theory:
 \begin{equation}
    \begin{tikzcd}
    K_0(\mathbb{K}(\mathcal{F})) \ar{r}{j_{\ast}} & K_0(\Toe) \ar{r}{q_{\ast}} & K_0 \ar{d}{\exp}(\Pim) \\
    K_1(\Pim) \ar{u}{\partial} & K_1 \ar{l}(\Toe) \ar{l}{q_{\ast}} & K_1(\mathbb{K}(\mathcal{F})) \ar{l}{j_{\ast}},
    \end{tikzcd} 
\end{equation}
Observe that the algebra of compact operators is Morita equivalent to the complex numbers via the Fock space $\mathcal{F}$ of the subproduct system. We shall denote the $KK$-class of the Fock space and its dual by $[\mathcal{F}]\in KK(\mathbb{K}(\mathcal{F}), \C)$ and $[\mathcal{F}^*]\in KK(\C, \mathbb{K}(\mathcal{F}))$, respectively.
\begin{prop}
Let $P_1$ and $P_2$ be Temperley--Lieb polynomials. Denote by $\mathcal{H} = \mathcal{H}_{P_1} \star \mathcal{H}_{P_2}$ the free product of their subproduct systems. The following identity 
    \begin{align*}
        [j] \otimes_\Toe [\psi_+, \psi_-] = [\mathcal{F}_{\mathcal{H}}] \otimes_\C \left(\textbf{1}_\C - [H_1] + [H_2^\perp]\right)
    \end{align*}
  holds  in $KK(\mathbb{K}, \C)$.
\end{prop}
\begin{proof}
    Since $H_2^\perp$ is $2$-dimensional, it is sufficient to show that 
    \begin{align*}
        [j] \otimes_\Toe [\psi_+, \psi_-] = 3\cdot [\mathcal{F}_\mathcal{H}] - [\mathcal{F}_\mathcal{H}]\otimes_\C [H_1].
    \end{align*}
 The proof is a simple adaptation of that of \cite[Proposition 7.1]{AK21}.
\end{proof}
The above proposition, combined with $KK$-equivalence proven in Theorem \ref{KK-equi-formula}, yields
\begin{equation}
\xymatrix{0 \ar[r] & K_1(\Pim_{\langle P_1,P_2 \rangle}) \ar[rr]^{([\mathcal{F}_\mathcal{H}]\otimes_\mathbb{K} \cdot )\circ \partial} & & K_0(\C) \ar[rrr]^{(\textbf{1}_\C - [H_1] + [H_2^\perp])} & & & K_0(\C) \ar[r]^-{(q\circ i)_{*}} & K_1(\Pim_{\langle P_1,P_2 \rangle}) \ar[r] & 0}, \label{gysin-sequence}
 \end{equation}
where $q\circ i: \C \to \Pim_{\langle P_1,P_2 \rangle}$

\begin{cor}\label{cor:KtheoryFP}
Let $P_1$ and $P_2$ be Temperley--Lieb polynomials in $d_1$ and $d_2$ variables, respectively. Then
    \[K_0(\Pim_{\langle P_1 , P_2 \rangle}) \cong \mathbb{Z}/(d_1+d_2-1)\mathbb{Z}, \qquad K_1(\Pim_{\langle P_1,P_2 \rangle}) \cong 
       \lbrace 0 \rbrace.\]
\end{cor} 
\noindent Note that the above $K$-groups are unchanged when one swaps $d_1$ with $d_2$.

By induction, we can extend this result to the case of finitely many Temperley--Lieb polynomials $P_1, \dots, P_r$, when the Fock space $\mathcal{F}_\mathcal{H}$ is the free product of $\mathcal{F}_{\mathcal{H}_{i}}, i = 1, \dots, r$, where $\mathcal{H}_{i}$ denotes the subproduct system associated with Temperley--Lieb polynomial $P_i$, and the Toeplitz algebra $\Toe_{\mathcal{H}}$ is $\ast$-isomorphic to the reduced free product of $\Toe_{\mathcal{H}_i}, i = 1, \dots, r$.

\begin{cor}\label{cor:KtheoryFPmany}
For $i=1, \dots, r$, let $P_i$ be a Temperley--Lieb polynomial in $d_i$ variables. Then 
    \[K_0(\Pim_{\langle P_1 , \dots, P_r \rangle}) \cong \mathbb{Z}/(\sum_{i = 1}^r d_i - 1)\mathbb{Z}, \qquad K_1(\Pim_{\langle P_1, \dots, P_r \rangle \rangle}) \cong 
       \lbrace 0 \rbrace.\]
\end{cor} 

\section{Subproduct systems from quantum group corepresentations} \label{sec:SUq2}
Our interest in the representation theory of $SU(2)$ and of its quantum counterpart, Woronowicz's $SU_q(2)$, stems from their importance in various fields within mathematical physics, where they play a crucial role both in the study of symmetries and in quantum mechanics. 

In \cite{AK21}, the authors gave a recipe for constructing a subproduct system of finite-dimensional Hilbert spaces starting from a finite-dimensional representation of the compact group $SU(2)$ on a Hilbert space $V$. In their construction, the main ingredient in their construction was the so-called determinant of the representation, a subspace of the vector space $V \otimes V$. We will provide here an alternative and more compact definition for that notion. The authors would like to thank Marcel de Jeu for pointing this out to us.

\begin{definition}
Let $(\rho, V)$ be a finite-dimensional unitary representation of the group $SU(2)$. Define the \emph{determinant} of the representation $(\rho, V)$ as the isotypical component of the trivial representation in $(\rho \otimes \rho, V\otimes V)$.     
\end{definition}

This definition can be dualised to the case of a corepresentation of the Hopf $C^*$-algebra $SU_q(2)$, and more generally, to the setting of a corepresentation of a rank-two compact quantum group.  We assume the reader to be familiar with the relevant notions from the theory of quantum group corepresentations \cite{NeTu13}, in particular with the tensor product of two corepresentations.

\begin{definition}
\label{def:det}
Let $\rho: V \to V \otimes C(SU_q(2))$ be a right corepresentation of the quantum group $SU_q(2)$ on $V$. We define the \emph{determinant} of $\rho$ as the isotypical component of the trivial corepresentation in the diagonal corepresentation $\rho \otimes \rho$ on the tensor product $H \otimes H$: 
\[
\mathrm{det}(\rho,V) = \lbrace \xi \in V \otimes V \mid 
\big( \rho \otimes \rho \big)(\xi) = \xi \otimes 1  \rbrace.
\]
\end{definition}
Note that since the determinant is a subspace of $H \otimes H$, taking its orthogonal complement gives a quadratic subproduct system of Hilbert spaces. 
\begin{example}
Recall that the fundamental corepresentation $\rho_1: \mathbb{C}^2 \to \mathbb{C}^2 \otimes C(SU_q(2))$, i.e. the irreducible corepresentation of $SU_q(2)$ with highest weight $1$, has matrix coefficients
\[\begin{pmatrix}
a & -qc^* \\
c & a^* 
\end{pmatrix}.\]
Let us consider the standard basis of $\mathbb{C}^2$. 
%Then the vectors $e_1$ and $e_2$ transform as:
%\[ e_1 \mapsto e_1 \otimes a + e_2 \otimes c, \qquad  e_2 \mapsto -q e_1 \otimes c^* + e_2 \otimes a^*.\]
It is easy to check that the determinant is spanned by the Temperley--Lieb vector
\begin{equation}
\label{eq:detsuq}
q^{-1/2} e_1 \otimes e_2 - q^{1/2} e_2 \otimes e_1 .\end{equation}
This follows from the commutation relations of $SU_q(2)$, in particular 
\[ac^*=qc^*a, \quad a^*a+c^*c=1 = aa^* + q^2 cc^*,\]
together with the fact that $c$ is normal.

Note that $\det(\rho_1,\mathbb{C}^2)$ is nothing but the $q$-antisymmetric subspace of $\mathbb{C}^2 \otimes \mathbb{C}^2$ defined in \cite{PuWo89} using the braiding $\sigma_q$ given by
\[\begin{cases}
\sigma_q(e_i \otimes e_i) = e_i \otimes e_i, & i=1,2 \\
\sigma_q(e_1 \otimes e_2) = q e_2 \otimes e_1, & \\
\sigma_q(e_2 \otimes e_1) = q e_1 \otimes e_2 + (1-q^2) e_2 \otimes e_1.
\end{cases}\]
Its orthogonal complement is the so-called $q$-symmetric tensor product.
\end{example}

It is a well-known fact that the group $SU(2)$ and its quantum analogue $SU_q(2)$ have the same representation category, and hence the same fusion rules.

\begin{theorem} {\cite[Theorem~5.11]{SW}}
    Let $V_n, V_m$ be the irreducible corepresentation of $SU_q(2)$ with highest weights $n$ and $m$, respectively. Then the tensor product of corepresentations $V_n \otimes V_m$ decomposes as
    \begin{align*}
        V_n \otimes V_m \cong V_{|n - m|} \oplus V_{|n - m| + 1} \oplus \dots \oplus V_{n + m}.
    \end{align*}
\end{theorem}
To characterize the determinant of $V_n$, we apply the Clebsch--Gordan formula {\cite[Equation~(54)]{KlKo}}, we have 
\begin{align} \label{eq:q-det}
    \text{det}(\rho_n,\mathbb{C}^{n + 1}) = \text{span}_\C \left\lbrace \sum_{i = 1}^{n+1} (-1)^i \cdot \left( \frac{q^{-n + i}}{[n + 1]_q} \right)^{1/2} \cdot e_i \otimes e_{n + 2 - i}, \right \rbrace
\end{align}
where $\{e_1, e_2, \dots, e_{n+1} \}$ is an orthonormal basis of $V_n \cong \C^{n+1}$, see also cf.\cite{GeW}. The vector above is Temperley--Lieb {\cite[Lemma~1.4]{HaNe21}}, with corresponding Temperley--Lieb polynomial.
\begin{equation}
\label{eq:TLSUq}
P(X_1, \dots, X_n) = \sum_{i = 1}^{n+1} (-1)^i \cdot \left(\frac{q^{-n + i}}{[n + 1]_q} \right)^{1/2}\cdot X_i X_{n + 2 - i}. \end{equation}

Therefore, the $SU_q(2)$-subproduct system is a Temperley--Lieb subproduct system.

\begin{lemma} \label{dim:su_q}
    Let $\rho$ be a finite-dimensional corepresentation of $SU_q(2)$, then the determinant has dimension equal to the sum of the squares of the multiplicities of its irreducible components.
\end{lemma}
\begin{proof}
    Let $\textbf{1}$ denote the trivial corepresentation. By Definition~\ref{def:det}.  $\det(\rho)$ is the isotypical component of $\textbf{1}$ in $(\rho \otimes \rho, H \otimes H)$, and it is thus determined by the intertwiner space $\text{Hom}(\rho \otimes \rho, \textbf{1})$. 
    
    Let $\rho$ be a finite-dimensional reducible corepresentation of $SU_q(2)$. Then $\rho$ decomposes into the direct sum of irreducible corepresentations $\rho_n$ of highest weight $n$, with multiplicity $k_n$, i.e., $\rho = \bigoplus_{n = 0}^{\infty} \rho_n^{\oplus k_n}$ with finitely many $k_n$'s non-zero. We compute
\[
    \text{Hom}(\rho \otimes \rho, \textbf{1}) \cong \text{Hom}(\rho, \rho^*) \cong \text{Hom}\left(\bigoplus_n \rho_{n}^{\oplus k_n}, \bigoplus_n (\rho_{n}^{\oplus k_n})^*\right).
\]
Given that $\rho_{n}$ is not equivalent to $\rho_{m}$ for $n \neq m$, and irreducible corepresentations of $SU_q(2)$ are self-dual, we obtain
\[    \text{Hom}(\rho \otimes \rho, \textbf{1}) \cong \bigoplus_n \text{Hom}\left(\rho_{n}^{\oplus k_n}, (\rho_{n}^{\oplus k_n})^*\right) 
    \cong \bigoplus_n \text{Hom}(\rho_{n}^{\oplus k_n}\otimes \rho_{n}^{\oplus k_n}, \textbf{1}).\] 
We deduce that
    \begin{align*}
        \dim(\text{Hom}(\rho \otimes \rho, \textbf{1})) &= \sum_{n = 0}^\infty \dim(\text{Hom}(\rho_{n}^{\oplus k_n}\otimes \rho_{n}^{\oplus k_n}, \textbf{1})) \\
        &= \sum_{n = 0}^\infty \dim(\text{Hom}((\rho_{n}\otimes \rho_{n})^{\oplus k_n^2}, \textbf{1})) = \sum_{n = 0}^\infty k_n^2,
    \end{align*}
    which completes the proof.
\end{proof}

\subsection{The subproduct system of a multiplicity-free corepresentation}
\begin{theorem}
\label{thm:multfree}
    Let $(\rho, H)$ be a finite-dimensional multiplicity-free corepresentation of $SU_q(2)$. The $SU_q(2)$-subproduct system of $\rho$ is isomorphic to the free product of the $SU_q(2)$-subproduct systems of its irreducible components. Correspondingly, the Toeplitz algebra $\Toe_{\mathcal{H}}$ is the reduced free product of the Toeplitz algebras of the subproduct systems of its irreducible components.
\end{theorem}
\begin{proof}
    To establish the result, it suffices to show that the determinant of a multiplicity-free unitary representation is spanned by Temperley--Lieb vectors. 
    %Indeed, if the Temperley spans the determinant--Lieb polynomials, the corresponding $SU_q(2)$-subproduct system is a quadratic subproduct system whose quadratic relations are precisely given by the Temperley--Lieb polynomials. By definition, the free product of quadratic subproduct systems is determined by the relations of its constituent subproduct systems. Therefore, the $SU_q(2)$-subproduct system of a finite-dimensional multiplicity-free corepresentation is the free product of Temperley--Lieb subproduct systems.
By definition, $\det (\rho)$ is the isotypical component of the trivial corepresentation within $(\rho \otimes \rho, H \otimes H)$ and is thus determined by the intertwiner space $\text{Hom}(\rho \otimes \rho, \textbf{1})$. Since $\rho$ is multiplicity-free, a similar argument to that in Lemma \ref{dim:su_q} gives
\begin{align*}
     \text{Hom}(\rho \otimes \rho, \textbf{1}) \cong \oplus_i \text{Hom}(\rho_{n_i}\otimes \rho_{n_i}, \textbf{1}).
\end{align*}
  Consequently, the determinant of the representation $\rho$ decomposes as the direct sum of the determinants of its irreducible components. 

    Let us decompose $\rho$ into its irreducible components, i.e. $\rho = \rho_{n_1} \oplus \dots \oplus \rho_{n_r}$, where $\rho_{n_i}$ denotes the irreducible corepresentation with highest weight $n_i$ and $\rho_{n_k} \not\cong \rho_{n_l}$ for $k \neq l$. Denote by $\mathcal{H}^i$ the $SU_q(2)$ subproduct system associated with $\rho_{n_i}$ and by $\mathcal{H}$ the $SU_q(2)$ subproduct system associated with $\rho$.
  For each irreducible component $\rho_{n_i}$, the determinant is one-dimensional and spanned by the vector in \ref{eq:q-det}. Therefore, the determinant of $\rho$ is spanned by a union of independent Temperley--Lieb vectors.
\end{proof}
Since the Toeplitz algebra associated with an irreducible representation of $SU_q(2)$ is nuclear \cite[Corollary~3.3]{HaNe21}, combining Theorem~\ref{thm:multfree} with Theorem \ref{thm:fp_kk-eq}, we deduce the following:
\begin{cor}\label{cor:mf_nuc}
    Let $(\rho, H)$ be a finite-dimensional multiplicity-free corepresentation of $SU_q(2)$, and $\mathcal{H}$ be the associated $SU_q(2)$-subproduct system.
Then the Toeplitz algebra $\Toe_{\mathcal{H}}$ is nuclear. 
\end{cor}

From the fact that the subproduct system of irreducible $SU_q(2)$-representation is Temperley--Lieb, we may view the exact sequence \eqref{gysin-sequence} as a noncommutative Gysin sequence. To this end, we define the Euler class of the representation to be
\[ \chi\left(\oplus_{i = 1}^r \rho_{n_i}\right) := \textbf{1}_\C - [H_1(\oplus_{i = 1}^r \rho_{n_i})] + [\det(\oplus_{i = 1}^r \rho_{n_i})] \in KK(\mathbb{C},\mathbb{C}).\]
\begin{theorem} \label{K-theory-CP}
    We have an exact sequence of groups:
    \begin{equation}
    \begin{tikzcd}
  0 \ar{r}&  K_1(\Pim) \ar{rr}{([\mathcal{F}]\otimes_\mathbb{K} \cdot )\circ \partial} & & K_0(\C) \ar{rrr}{[\mathcal{F} ] \otimes_\C \chi(\oplus_{i = 1}^r \rho_{n_i})} & & &  K_0(\C) \ar{rr}{(q\circ i)_*} & & K_1(\Pim) \ar{r}& 0
    \end{tikzcd}, \label{gysin-sequence-concrete}
\end{equation}
%where $\mathcal{F}$ denotes the Fock space of $\mathcal{H}$.
Therefore, we have  
\begin{align*}
    &K_1(\Pim) \cong \ker \left( \textbf{1}_{\C} - [E_1(\oplus_{i = 1}^r \rho_i)] + [\det(\oplus_{i = 1}^r \rho_i)])\right),\\
    &K_0(\Pim) \cong \coker\left( \textbf{1}_\C - [E_1(\oplus_{i = 1}^r \rho_i)] + [\det(\oplus_{i = 1}^r \rho_i)]) \right).
\end{align*} 
More precisely, for $\rho \cong \oplus_{i = 1}^r \rho_{n_i}$, the $K$-theory groups of its Cuntz--Pimsner algebra are
\begin{align*}
    K_0(\Pim) \cong \mathbb{Z}/(\sum_{i = 1}^r n_i - 1)\cdot \mathbb{Z} \quad \quad K_1(\Pim) \cong \begin{cases}
        \mathbb{Z} \quad r = 1, n_1 = 1\\
        0 \quad \text{otherwise}.
    \end{cases}
\end{align*}
\end{theorem}
\begin{remark}
    The above theorem extends \cite[Corollary~7.3]{AK21} beyond the irreducible case.
\end{remark}

\subsection{Dealing with multiplicities}
Let $\rho$ be an isotypical corepresentation of $SU_q(2)$ with highest weight $n$ and multiplicity $t$, i.e. $\rho = \rho_n^{\oplus t}$, then by Lemma \ref{dim:su_q}, we have $\dim(\det(\rho)) = t^2$. Indeed, an explicit basis for $\det(\rho)$ is the following:

\[ \text{det}(\rho) = \left\lbrace \sum_{i = 1}^{n+1} (-1)^iq^{i/2} e^k_i\otimes e^l_{n + 1 - i} : k, l = 1, 2, \dots, t  \right\rbrace, \]
where we denote \[ e_i^k = 0 \oplus  \dots 0 \oplus \underbrace{e_i}_{k \text{th}} \oplus 0 \oplus \dots \oplus 0,\]and the common divisor of the coefficient $[n + 1]_q^{- 1/2}$ is omitted. 

\begin{remark}
\label{rem:genIso}
The subproduct system of an isotypical representation is a quadratic subproduct system with few relations. Indeed, let $\dim(\rho)=m$ then we have $mt$ generators and $t^2$ relations, and it is easy to see that $t^2 \leq (mt)^2/4$ precisely when $m \geq 2$.
\end{remark}
Let $n \in \N$. By $\rho_n$ we mean the irreducible corepresentation of $SU_q(2)$ of highest weight $n$. For simplicity, we denote the associated $SU_q(2)$-subproduct system $\mathcal{H}:= \{ H_m \}_{m \in \N_{0}}$, where $H_1$ is the representation space. Moreover, we denote the $SU_q(2)$-subproduct system associated with the representation $\rho = \rho_n^{\oplus t}$ by $\{ H^t_m : m \in \N \}$. For any $1 \leq k \leq t$, we define $\sigma^1_k:  H_1 \to {H}^{t}_{1}$ as the $SU_q(2)$-equivariant linear maps given on the basis vectors by $\sigma^1_k(e_i) = e^k_i$. 

By definition, ${H}_1^t \cong H_1^{\oplus t}$ through the $SU_q(2)$-equivariant isomorphism given by
\[ {H}_{1}^{t} \cong \sigma^1_1(H_1) \oplus \sigma^1_2(H_1) \oplus \dots \oplus \sigma^1_t(H_1). \]

The vector space $H^{t}_m$ can be described in a similar way: 
%
% More generally, we have the following proposition.
\begin{prop}
    Let $\rho_n$ denote the irreducible corepresentation of highest weight $n$ and $\mathcal{H}$ be the corresponding $SU_q(2)$-subproduct system. Let $\mathcal{H}^t$ be the subproduct system of the corepresentation $\rho_n^{\oplus t}$. There is a unitary isomorphism:
    \begin{equation}
        \label{eq:iso_isotypical}
      % \bigoplus_{k_1, k_2, \dots, k_m = 1}^t \sigma^1_{k_1}\otimes \sigma^1_{k_2} \otimes \dots \otimes \sigma_{k_m}^1 :  
      H_m^{\oplus t^m} \simeq H^{t}_m     
    \end{equation}
    \end{prop}
\begin{proof}
    We will show that the isomorphism is implemented by the map 
    \[ \bigoplus_{k_1, k_2, \dots, k_m = 1}^t \sigma^1_{k_1}\otimes \sigma^1_{k_2} \otimes \dots \otimes \sigma_{k_m}^1 :  
      H_m^{\oplus t^m} \to H^{t}_m.\]
      We prove this by induction. The statement is true for $n=1$. For $n=2$, recall the definition of ${H}^{t}_2$ as the orthogonal complement of the determinant in $H^t_1\otimes H^t_1 $. Observe that 
 \[ H^t_1\otimes H^t_1 \cong \bigoplus_{k_1, k_2 = 1}^t \sigma^1_{k_1}(H_1) \otimes \sigma^1_{k_2}(H_1). \]
 Moreover, we have that \[D:= \det (\rho_n^{\oplus t})\cong \bigoplus_{k_1, k_2 = 1}^t \sigma^1_{k_1} \otimes \sigma^1_{k_2} (\det(\rho_n)).\]Therefore, 
 \begin{align*}
  H^t_2 &= D^\perp \cong \bigoplus_{k_1, k_2 = 1}^t \sigma^1_{k_1}(H_1) \otimes \sigma^1_{k_2}(H_1) \ominus \bigoplus_{k_1, k_2 = 1}^t \sigma^1_{k_1}
     \otimes \sigma^1_{k_2} (\text{det}_q(\rho_n)) \\
    % H^t_2 &= H^t_1\otimes H^t_1 \ominus D \\
     %&\cong \bigoplus_{k_1, k_2 = 1}^t \sigma^1_{k_1}(H_1) \otimes \sigma^1_{k_2}(H_1) \ominus \bigoplus_{k_1, k_2 = 1}^t \sigma^1_{k_1} \otimes \sigma^1_{k_2} (\text{det}_q(\rho_n)) \\
 &\cong \bigoplus_{k_1, k_2 = 1}^t  \sigma^1_{k_1} \otimes \sigma^1_{k_2}(H_1 \otimes H_1 \ominus \det(\rho_n)) \cong \bigoplus_{k_1, k_2 = 1}^t  \sigma^1_{k_1} \otimes \sigma^1_{k_2}(H_2),
 \end{align*}
which proves the claim for $m=2$.

Using the recursive formula in Remark \ref{rem:intersection}, we obtain
    \begin{align*}
      &{H}^t_{m+1} \cong  H^t_1\otimes H^t_{m} \cap H^t_{m} \otimes H^t_1 \\
       & \quad \cong \left(\sigma^1(H_1)^{\oplus t} \otimes (\bigoplus_{k_1, \dots, k_m = 1}^t \sigma^1_{k_1} \otimes \dots \otimes \sigma_{k_m}^1 (H_m))\right) \cap \left((\bigoplus_{k_1, \dots, k_m = 1}^t \sigma^1_{k_1} \otimes \dots \otimes \sigma_{k_m}^1 (H_m)) \otimes \sigma^1(H_1)^{\oplus t}\right) \\
       % &\cong \bigoplus_{k, k_1, \dots, k_m = 1}^t \sigma^1_k\otimes \sigma^1_{k_1} \otimes \dots \otimes \sigma_{k_m}^1 (E_1 \otimes E_m) \\
        %&\cap \oplus_{k, k_1, \dots, k_m = 1}^t \sigma^1_{k_1}\otimes \sigma^1_{k_2} \otimes \dots \otimes \sigma_{k_m}^1 (E_m \otimes E_1) \\
        &\quad \cong \bigoplus_{k, k_1, \dots, k_m = 1}^t \sigma^1_k\otimes \sigma^1_{k_1} \otimes \dots \otimes \sigma_{k_m}^1 (H_1 \otimes H_m \cap H_m \otimes H_1) \\
        & \quad \cong \bigoplus_{k, k_1, \dots, k_m = 1}^t \sigma^1_k\otimes \sigma^1_{k_1} \otimes \dots \otimes \sigma_{k_m}^1 (H_{m+1}).
    \end{align*}
    which concludes the proof.
\end{proof}

\begin{cor}
Let $\rho_n$ be the irreducible $SU_q(2)$ corepresentation with highest weight $n$, and let $h_n(z)$ be the Hilbert series of the associated $SU_q(2)$ subproduct system. The Hilbert series of the subproduct system of the isotypical corepresentation $\rho_{n}^{\oplus t}$ satisfies
\begin{equation}
\label{eq:HSmult}
h_{tn}(z) = (1-t(n+1)z + t^2z^2)^{-1} = h_{n} (tz) \end{equation}
\end{cor}
\begin{proof}
The proof follows from the corresponding claim for dimension sequences: let $d^{(n)}$ be the dimension sequence of the subproduct system of the irreducible corepresentations $\rho_n$. Then the subproduct system of the isotypical corepresentation $\rho_n^{\oplus t}$ is given by
$d^{(n,t)}_m:= d^{(n)}_m t^{m}. $
Our claim then follows from the definition of Hilbert series.
\end{proof}
Combining this result with Remark~\ref{rem:genIso}, we obtain the following:
\begin{cor}
Let $\rho_n$ be the irreducible $SU_q(2)$-corepresentation of highest weight $n$. The subproduct system of the corepresentation $\rho_{n}^{\oplus t}$  is a generic quadratic subproduct system in $t(n+1)$ generators and $t^2$ relations.
\end{cor}
\begin{example}
    Let $\rho_1$ be the fundamental corepresentation on  $V_1 \simeq \C^2$ with orthonormal basis $\{ e_1, e_2 \}$.  Then the determinant is given by \eqref{eq:detsuq}.
    
    Consider the isotypical corepresentation $\rho_1^{\oplus 2}$ on $H_1^t$ with basis $\{ e_1^1, e^1_2, e^2_1, e^2_2 \}$. Then we have $\det(\rho_1^{\oplus 2}) \cong \det(\rho_1)^{\oplus 4}$, which is spanned by \[\{ q^{1/2}\cdot e^k_1\otimes e^l_2 - q^{-1/2} \cdot e^k_2 \otimes e^l_1: k, l = 1, 2 \}.\] 
The space $H_2^t \cong H_2^{\oplus 4}$ is spanned by
\[ \{ e^k_1 \otimes e_1^l, e^k_2 \otimes e^l_2, q^{-1/2} \cdot e^k_1\otimes e^l_2 + q^{1/2}\cdot e^k_2\otimes e^l_1 : k, l = 1, 2 \}  . \]

As discussed earlier, this construction gives a generic quadratic subproduct system with few relations, with Hilbert series 
\[ h(z)= (1-4z+4z^2)^{-1}.\] 
\end{example}

\subsection*{Outlook}
It is natural to wonder what operation in the algebraic world of associative algebras corresponds to the change of variable in the Hilbert series described in \eqref{eq:HSmult}, and to consider what the consequences of this operation are at the level of the Toeplitz algebras. 

Finally, it seems that K-theory computations only read the Hilbert series of a quadratic algebra and that Cuntz--Pimsner algebras of non-isomorphic subproduct systems with the same Hilbert series are KK-equivalent. We postpone the discussion of these and other related questions to future work.
%We shall now give a more explicit description of the $SU(2)$-equivariant unitary isomorphism
%\begin{align}
%    W^n_R: E_n\otimes E_1 \cong E_{n+1} \oplus E_{n-1}^{\oplus R}. \label{eq:fusion}
%\end{align}

\end{document}